\documentclass[preprint, 12pt]{elsarticle}
\usepackage[utf8]{inputenc}

\usepackage{todonotes}
\usepackage{caption}
\usepackage{fancyhdr}

\usepackage{hyperref}
\usepackage{amsfonts}
\usepackage{geometry}
\usepackage[ruled]{algorithm2e}

\usepackage{amsmath}

\usepackage{amssymb}   
\usepackage{amsthm}

\usepackage{pict2e}
\usepackage{keyval}
\usepackage{calc}
\usepackage{fp}
\usepackage{diagbox}
\usepackage{booktabs}
\usepackage{tabularx,colortbl, makecell, caption, multirow}
\usepackage{boldline} 

\newtheorem{theorem}{Theorem}[section]

\newtheorem{definition}[theorem]{Definition}

\usepackage{rotating}
\usepackage{graphicx}
\usepackage{bm}
\usepackage{xcolor}
\usepackage{multirow}
\usepackage{nicematrix}

\usepackage{graphicx}
\graphicspath{{images/}}
\usepackage{float}
\usepackage{longtable}
\usepackage{adjustbox}

\usepackage{colortbl} 

\graphicspath{{./figures/}}

\usepackage[titletoc,toc,title]{appendix}
\usepackage{comment}
\usepackage{graphics}
\usepackage{graphicx}
\usepackage{adjustbox}
\usepackage{tikz}
\usepackage{caption}
\usepackage{psfrag}
\usepackage{soul}
\usepackage{mwe}
\usepackage{pgfplots}
\pgfplotsset{compat=1.17} 

\graphicspath{{./figures/}}

\usepackage{comment}
\usepackage{graphics}
\usepackage{graphicx}
\usepackage{adjustbox}
\usepackage{tikz}
\usepgfplotslibrary{dateplot}
\usepgfplotslibrary{groupplots}
\usepackage{caption}
\usepackage{subcaption}
\usepackage{psfrag}
\usepackage{soul}
\usepackage{mwe}
\usepgfplotslibrary{fillbetween}

\usetikzlibrary{er,positioning, calc, shapes, arrows, chains, fit}

\definecolor{mygreen}{RGB}{80,176,50}

\usepackage{xcolor}

\makeatletter

\newcommand*{\@rowstyle}{}

\newcommand*{\rowstyle}[1]{
  \gdef\@rowstyle{#1}%
  \@rowstyle\ignorespaces%
}

\newcolumntype{=}{
  >{\gdef\@rowstyle{}}%
}

\newcolumntype{+}{
  >{\@rowstyle}%
}

\makeatother
\def\T{\mathcal{T}}
\def\B{\mathcal{B}}
\def\Tb{\mathcal{B}}
\def\Tbf{\mathcal{B}'}
\def\Tbl{\mathcal{B}''}

\def\Tl{\mathcal{L}}
\def\L{\mathcal{L}}

\def\D{\mathcal{D}}
\def\H{\mathcal{H}}
\def\S{\mathcal{S}}

\def\I{\mathcal{I}}

\def\F{\mathcal{F}}
\def\J{\mathcal{J}}
\def\U{\mathcal{U}}
\def\E{\mathcal{E}}

\newcommand{\R}{\mathbb{R}}

\newcommand{\ds}{\displaystyle}

\definecolor{dred}{rgb}{0.55, 0.0, 0.0}
\definecolor{dblue}{rgb}{0.0, 0.2, 0.6}
\definecolor{dgreen}{RGB}{80,176,50}
\definecolor{mygreen}{RGB}{80,176,50}
\definecolor{myred}{RGB}{222,36,16}
\definecolor{myind}{RGB}{93,174,255}
\definecolor{mypink}{RGB}{219,129,255}
\definecolor{azzdark}{RGB}{42,72,180}
\definecolor{azzlight}{RGB}{53,146,203}
\definecolor{verdone}{RGB}{0,128,0}
\definecolor{indaco}{RGB}{68,137,206}
\definecolor{celestino}{RGB}{222,247,254}
\definecolor{verdeacceso}{RGB}{71,185,109}
\definecolor{pred}{RGB}{182,7,40}
\definecolor{pblue}{RGB}{58,76,189}
\definecolor{orange}{rgb}{0.8, 0.47, 0.13}
\definecolor{green}{rgb}{0.0, 0.62, 0.42}

\newcommand{\dblue}[1]{\textcolor{dblue}{#1}}

\newcommand{\verdone}[1]{\textcolor{verdone}{#1}}

\newcommand{\mnote}{\textcolor{black}}
\newcommand{\gnote}{\textcolor{black}}
\newcommand{\lnote}{\textcolor{black}}
\newcommand{\lcomment}{\textcolor{black}}

\renewcommand{\thanks}[1]{}
\renewcommand{\footnote}[1]{}

\begin{document}


\begin{frontmatter}

\title{
Unboxing Tree Ensembles for interpretability: \\
a hierarchical visualization tool and\\  a multivariate optimal re-built tree 
}

\author[1]{Giulia Di Teodoro}
\author[1]{Marta Monaci}
\author[1]{Laura Palagi}

\address[1]{Department of Computer, Control and Management Engineering Antonio Ruberti (DIAG), \\Sapienza University of Rome, 00185, Rome, Italy\\
          
\{giulia.diteodoro;marta.monaci;laura.palagi\}@uniroma1.it
            }

\begin{abstract}

{The interpretability of models has become a crucial issue in Machine Learning because of algorithmic decisions' growing impact on real-world applications. Tree ensemble methods, such as Random Forests or XgBoost, are powerful learning tools for classification tasks. However, while combining multiple trees may provide higher prediction quality than a single one, it sacrifices the interpretability property resulting in "black-box" models. In light of this, we aim to develop an interpretable representation of a tree-ensemble model that can provide valuable insights into its behavior. First, given a target tree-ensemble model, we develop a hierarchical visualization tool based on a heatmap representation of the forest's feature use, considering the frequency of a feature and the level at which it is selected as an indicator of importance.
Next, we propose a mixed-integer linear {programming (MILP)} formulation for constructing a single optimal multivariate tree that accurately mimics the target model predictions. The goal is to provide an interpretable surrogate model based on oblique hyperplane splits, which uses only the most relevant features
according to the defined forest's importance indicators. 
{The MILP model includes a penalty on feature selection based on their frequency in the forest to further induce sparsity of the splits.} 
The natural formulation has been strengthened to improve the computational performance of {mixed-integer} software. Computational experience is carried out on benchmark datasets from the UCI repository using a state-of-the-art off-the-shelf solver. Results show that the proposed model is effective in yielding a shallow interpretable tree approximating the tree-ensemble decision function.} 

\end{abstract}


\end{frontmatter}

\parindent 0cm


\pagestyle{fancy}
\lhead{}
\renewcommand{\headrulewidth}{0pt}
\rhead{\textcolor[gray]{0.3}{This version has been published in EJCO: \href{https://doi.org/10.1016/j.ejco.2024.100084}{https://doi.org/10.1016/j.ejco.2024.100084}}}
\thispagestyle{fancy}  

\textcopyright{ 2024. Licensed under the Creative Commons  \href{https://creativecommons.org/licenses/by-nc-nd/4.0/}{CC-BY-NC-ND 4.0}.} 

\section{Introduction}

When building Machine Learning (ML) models in supervised learning, it is becoming more and more important to balance accuracy and interpretability.
There is a growing number of sensitive domains where a detailed understanding of the model and the outputs is as important as the accuracy in prediction  \cite{burkart2021survey}.
Following \cite{doshi2017towards}, interpretability is defined as "the ability  to explain or to present the decision in understandable
terms to a human". 
Generally, models with high accuracy are more complex and less interpretable. This trade-off is evident when comparing Decision trees (DT) and Tree ensembles (TE) (e.g., Random Forests, XGBoost).
Decision Trees (DT) (see \cite{Breiman1984CART,quinlan:induction,Quinlan1993C45}) have high interpretability since their construction process is simple, intuitive,
 and can be easily visualized. They are white-box models easily explained by Boolean logic; indeed, the prediction made by the tree  is a conjunction of predicates.
The disadvantages are that they often overfit, do not have good out-of-sample predictive capabilities, can grow exponentially, and have high variability. 
On the other hand,   Tree Ensemble (TE) models such as Random Forests \cite{breiman2001random}, and XGboost \cite{Friedman2001}\cite{XGBoostChen} constitute one of the most widely used techniques for regression and classification tasks.
By aggregating many decision trees, tree ensemble techniques substantially increase predictive capabilities.
Indeed, including the randomness element and the ensemble procedure used to aggregate the individual tree predictions reduces the model variance. Thus, TEs can achieve high levels of accuracy  at the expense of lower interpretability. Indeed, a TE is a black-box model, not interpretable through parameters or its functional form.

TE models are used in many fields, such as medical or financial, where opaque and redundant decisions could be highly harmful, and an easily understood interpretation of the model's predictions could significantly impact the final decision.
Hence the ability to understand the interactions between predictors and responses used in a TE model is an important issue
when undertaking a decision-making process.

In  \cite{guidotti2018survey}, a list of desiderata of a surrogate predictive model is reported as follows. 
\begin{enumerate}
    \item Interpretability: to what extent the model and/or its predictions are human-understandable. We adopt as a measure of interpretability the complexity of the
predictive model in terms of the model size (depth or the number of leaves) and the sparsity of the predictors used to construct the decision rules.
\item Accuracy: to which extent the model accurately predicts out-of-sample instances. 
\item Fidelity: to what extent can the model accurately imitate a black-box predictor. 
\end{enumerate}
Measures of Accuracy and Fidelity are based on standard KPIs such as error rate, and F1-score having as target values the ground truth and the predicted values, respectively.



The paper is organized as follows. In Section \ref{sec:sota}, we review state of the art in Interpretative models for Tree Ensemble, including the main feature's VI indicators. In Section \ref{sec:contribution}, we focus on our contribution. In Section \ref{sec:basic}, we state the problem and the main basic definitions and tools used. Section \ref{sec:vite} is devoted to the description of the visualization toolbox, and Section \ref{sec:model} defines the model for the optimal re-built tree. The formulation is strengthened in Section \ref{sec:strenghen}.
Finally Section \ref{sec:numerical} reports numerical results on a standard benchmark of test problems from the UCI collection.

   

\section{State of the art on Interpretative models for Tree Ensemble}\label{sec:sota}
Several tools and methods in the literature are used to visualize and interpret tree ensembles, \lnote{ as the many trees in the forest make it difficult to understand the decision path leading to predictions.}

In this section, we aim to review only some of the possible approaches, and we refer to the recent survey 
\cite{ARIA2021100094} for a more detailed comparison. According to \cite{ARIA2021100094}, 
 the approaches can be mainly divided into: 
\begin{itemize}
    \item \textit{Internal processing}, which includes methods that aim to  provide a global overview of the model through measures helpful in interpreting the results obtained;
\item \textit{Post-Hoc approaches}, which aim to identify a relationship structure among response and predictors and include the construction of a single surrogate model that approximates the original TE {prediction function}.
\end{itemize}

\paragraph{Internal processing}
Among measures for interpreting the prediction results of a TE, the feature importance approach assigns each feature $j$  a score to indicate its impact on predicting the output $y$. 
These approaches produce a features ranking which is usually obtained by modifying in some way the value of each feature at the time and evaluating the impact on the quality of the tree ensemble.
Mean Decrease accuracy (MDA), Mean Decrease Impurity (MDI) \cite{breiman2001random,GENUER20102225,louppe2013understanding} and  Partial Dependence Plot (PDP) \cite{Friedman2001} provide examples of such measures.

A different approach to calculating the importance of variable $j$, presented in \cite{ishwaran2008random}, consists 
in randomly changing the assignment of the samples to the children nodes whenever the variable $j$ is used in the splitting rule of a node and  evaluating the change in the accuracy of the predictions. The more significant the change in the prediction, the more the importance of the variable $j$. 
{The underlying idea is that a random assignment of samples to nodes at shallower depths leads to more significant perturbations to the final assignment of samples to leaves and, thus, in the classification. 
}

In \cite{Ishwaran2012}, the Minimal Depth (MD) is introduced as an alternative measure of the importance based solely on the structure of the trees in the forest rather than the goodness of fit of the prediction. The features' importance is determined by the level at which they are used for splitting in the trees' nodes, where the more frequent use at shallow depths, the more important. The idea is that the higher the probability of a feature appearing in the first levels of the trees, the larger the impact on accuracy and also the importance accordingly to the noisy definition in \cite{ishwaran2008random}.
The Surrogate Minimal Depth, proposed in \cite{Seifert2019}, generalizes MD by using  surrogate variables, as defined by Breiman in \cite{Breiman1984CART}, to obtain the minimal depth.




Another way to provide additional information about the tree ensemble is the \emph{Proximity measure} of pair of samples, which is
a weighted average of the number of trees in the TE model in which
the samples end up in the same leaf. This concept was introduced first for Random Forest 
 \cite{breimancutler}, where weights can be all equal, and later extended to Gradient Boosted Trees \cite{Tan2020}.
The proximity measure takes values between 0 and 1. 
Associated with proximity, a similarity measure among samples can be defined as \emph{TE distance}  of pairs of samples. The TE distance is calculated as (1- TE proximity) and is a pseudo-metric \cite{Tan2020}.
The basic idea is that similar observations should be found more
frequently in the same terminal nodes than dissimilar observations.
Based on this measure, training samples can be visualized through a Multi-dimensional Scaling (MDS) plot so that users can intuitively observe data clusters and
outliers identified by the TE model.
Other visualization toolkits have been developed to graphically represent the measures shown above specifically for Random Forests:  \textit{randomForest} tool in R \cite{Liaw2002},  \textit{ggRandomForests} package in R \cite{ggRandomForest2015Ehrlinger}, and more recently  \textit{iForest} tool that allows interactive visualization of RF in Python  \cite{Zhao2018}.
These tools provide visualization of feature importance and partial dependence information, ranges of splitting values for each feature, similarities, and structure of decision paths, distribution of training data, and an interactive inspection of the model.

\paragraph{Post-Hoc approaches}
Another methodology used to aid the interpretability of TEs is represented by post-hoc approaches that aim to identify the relationship between the predicted output and the predictors. These include approaches that do \textit{size-reduction} of the forest, those that \textit{extract rules}, and those that aim for \textit{local explainability}. 

The model we propose in this paper falls into the size-reduction class. 
These approaches aim to reduce the size of TE while maintaining its predictive capability. This is an NP-hard problem \cite{Tamon2007}; in some cases, the smaller ensemble may even perform better \cite{ZHOU2002}. 
A huge branch of literature is available on these topics, and we refer to \cite{ARIA2021100094} and references therein for a recent review of the main approaches.

We are interested mainly in the stream of "born-again" trees, which involves building a single decision tree, called \emph{representer tree},  that mirrors the behavior of a pre-existing tree ensemble over all its feature space. Born-again trees were originally proposed in \cite{Breiman1996} as the problems of: "giving a probability distribution $\cal P$ on the space of input variables $x$,  find the tree that best represents $f(x)$" where  $f(x)$ is a given prediction.
In \cite{Breiman1996}, TEs are used to manufacture new pairs $(x^i,y^i)$, and a representer tree is constructed with a prediction accuracy on the manufactured samples close to the tree ensemble's accuracy rather than the ground truth. However, the tree can grow too much and might result in a not interpretable model. 
\lnote{From an optimization viewpoint, in \cite{vidal2020born}, a dynamic programming algorithm is proposed to build a born-again tree of minimal size that faithfully reproduces the behavior of the tree ensemble.}
\emph{Faithfulness} is defined as the capacity of the representer tree to reproduce precisely the decision function $f(x)$ of the TE on the whole space of the features and not only on the samples.
In both approaches \cite{Breiman1996,vidal2020born}, a univariate tree is constructed where the depth or number of leaves is controlled by optimizing exactly or approximately the complexity and by post-pruning to reduce the size. 
In \cite{vidal2020born}, it is proved that the problem of determining a faithful optimal decision tree of minimal size (where either depth or number of leaves, or any hierarchy of these two are used as complexity measures) is $\cal NP$-hard and that the depth of such an optimal tree is bounded above by the sum of the depth of the trees in the forest (the bound is tight). 
As the Problem is {\cal NP}-hard, the computational time of the dynamic programming algorithm proposed in \cite{vidal2020born} will eventually increase exponentially with the number of features, and a heuristic approach is proposed, which is guaranteed to be faithful but not necessarily minimal in size. A post-pruning phase is used on the set of six problems studied that seems not to affect the quality of the predictions but to significantly simplify the born-again trees.
Authors observed that born-again decision trees contain many inexpressive regions which do not contribute to effectively classifying samples.
The purpose of these regions and their contribution to the generalization capabilities of random forests is not clear yet.

\lcomment{Recently, in \cite{aria2023explainable}, an algorithm is proposed with the aim of representing an ensemble model using a dendrogram-like structure. The construction is based on the definition of a dissimilarity matrix, the complement of the proximity measure, to assess the discordance of observations with respect to the classifier. It constructs a tree-like structure for a clearer interpretation of the ensemble model.} 
\bigskip

\lcomment{\textit{Optimal classification trees}. It is well-known that learning an optimal binary decision tree is an \cal{NP}-complete problem \cite{Hyafil1976ConstructingOB}. 
For this reason, traditional approaches for building decision trees are based on sequential greedy heuristics, such as CART (Classification and Regression Trees) \cite{Breiman1984CART}. 
These approaches produce univariate decision trees where each splitting rule selects a single feature and an associate threshold to partition samples among the children nodes to minimize a local impurity function (e.g., the Gini index).
While these heuristic methods are computationally efficient, their greedy nature can lead to poor generalization performances.}

\lcomment{In recent years, thanks to the great improvement in mixed-integer programming (MIP), several works have been devoted to defining optimal classification trees (OCTs) using MIP methods. After the very first paper \cite{Bennett1992}, the seminal paper \cite{Bertsimas2017OptimalClassification} introduced
mixed-integer linear models for constructing optimal trees both with univariate splits and with multivariate ones.
This research paved the way for a variety of other formulations, which are reviewed in the comprehensive survey \cite{Morales2021MathematicalOptimizationOCRT}. Among recent papers focusing on multivariate OCTs,  in \cite{boutilier2022shattering}, a new MILP formulation and a new class of valid inequalities to
improve the optimization process are presented,  and in \cite{MarginOptimalClassificationTrees}, a new mixed-integer quadratic model based on Support Vector Machine is proposed for building optimal trees with maximum margin hyperplanes.}

\section{Our contribution}\label{sec:contribution}

Our contribution is twofold and follows the line of research aimed at creating an interpretable representation of a tree-ensemble model able to provide valuable insights into its behavior.
Following the classification proposed in \cite{ARIA2021100094}  of interpretative methods for Random Forests, we can divide our contribution into two parts interacting with each other.

In particular, given a Tree Ensemble (TE) model, we provide \begin{itemize}
    \item \texttt{VITE}: a hierarchical VIsualization tool for TE that depends solely on the structure of the tree and aims to visualize how features are used within the forest; 
    \item \texttt{\texttt{MIRET}}: a surrogate Multivariate Interpretable RE-built optimal Tree to gain interpretability on the relation between the input features and the TE outcomes.
\end{itemize}
 
The visualization tool \texttt{VITE} fits into the framework of the visualization tools developed in the literature to better understand the dynamics inside the TE. 
We drew our inspiration from the concept of \emph{Minimal Depth}, where the feature's importance is determined by the smallest depth at which it is used in a node's splitting rule for the first time. 
Indeed, since variable selection at splitting nodes is based on impurity indicators, those more frequently used at high levels produce the most significant impurity decrease.
However, we are interested in having a hierarchical view of the role of features in the trees composing the TE. Indeed, in a  tree, the \emph{Minimal Depth} (MD) of a feature $j$ is a nonnegative random variable 
taking values $\{0,\dots, D\}$, where $D$ is the depth of the tree. MD measures the distance from the root node to the root of the
closest maximal subtree having $j$ as a splitting feature. Essentially, it measures how far a sample moves down in the tree before encountering the first split on $j$.
In this way, the information on how often a feature is used in the tree, namely how many subtrees having $j$ as root nodes are present, is wholly lost. 
Indeed, it is possible for a feature to appear multiple times at lower levels in the tree after its initial appearance, which defines the Minimum Depth, or conversely, it may not appear again after its first appearance.
Following this idea,
we propose using the frequency of each feature's usage at each node or level in the trees of the TE as a measure of the feature's importance. This is based on the notion that features appearing more frequently at nodes closer to the root are likely to play a significant role in prediction outcomes.

We derive a tool for visualizing TE features' frequency on a single tree that gives{ a view at a glance of the hierarchical structure of the TE}.
Our tool is based on the heatmaps construction considering features' frequency. The TE visualization of \texttt{VITE} is intended to be an addition to the existing visualization tools for further facilitation of model interpretation in a graphic fashion.
\texttt{VITE} gives an immediate glimpse of the overall structure of the forest without the need to analyze individual decision pathways to extrapolate what might be the reason for different classifications.
In addition, it allows for seeing in a single view which are the features used most in the trees of the forest and the split point range of each feature, allowing one to reflect on the different ranges' effect on the predictions.

\par\medskip

The Multivariate Interpretable RE-built Tree \texttt{\texttt{MIRET}} fits in the framework of born-again trees, defining an \emph{optimal representer tree}.
A representer tree aims to replace a tree ensemble classifier with a newly constructed single tree that can reproduce, in some sense, the behavior of the TE. 
Differently from preceding approaches in \cite{Breiman1996,vidal2020born}, we propose constructing a tree with fixed depth, thus preventing it from growing excessively and becoming difficult to understand.
To this aim, we consider surrogate trees with fixed depth that uses oblique splits to partially regain the freedom lost by fixing complexity.

More in detail, based on the target TE, we present a mixed-integer linear programming (MILP) formulation for learning a Multivariate Interpretable RE-built tree
with the \emph{same maximum depth $D$} of TE, and that accounts for the information derived  from the target TE.

To this end, we extract knowledge from the target TE and inject information into our tree model. In particular:
\begin{itemize}
\item Following the TE voting procedure, which leads to the prediction, we maximize the fidelity to the TE. Namely, we minimize the misclassification of each sample with respect to its predicted class extracted from the ensemble model.
\item According to the measures used in the visualization tool, we detect each feature's usage frequency along the forest's trees.
 Through penalization in the objective function, we drive the selection of features based on their frequencies. Further, we hand to our model only a subset of the most frequent features for each tree level. In this way, we promote selecting the most representative features of the TE while further inducing the sparsity of the branching hyperplanes, i.e. sparsity of the features.
Indeed, sparsity is a core component of interpretability \cite{Rudin2022interpretable}, and having sparser decision rules allows the end user to identify better the key factors influencing the outcome.


\item {To encourage the partition of the feature space as performed by the TE}, we calculate the proximity of each sample pair in the TE. Based on this, we ensure that pairs of samples with proximity greater than a specified threshold are placed in the same final leaf of our tree, i.e., in the same final feature space partition.
\end{itemize}



In section \ref{sec:strenghen}, a strengthened MILP formulation is proposed to improve the computational performance of MILP algorithms.
This approach provides a simple yet effective way to interpret the predictions of a complex ensemble model, yielding a shallow interpretable tree able to give insights about the features that most affect the classification while approximating the RF decision function.

\section{Basic definition and preliminaries}\label{sec:basic}

We focus on binary classifiers, and we assume that we are given a training dataset $$\{(x^i,y^i) \in \mathbb{R}^{|\J|}\times\{-1,1\},\ i\in\I\},$$ where $\J$ is the index set of the features, and 
w.l.o.g. we normalize the feature values $x^i, i\in\I$ to lie inside the interval $[0,1]$.      

The training data is used to construct a predictor $f(x): \mathbb{R}^{|\J|}\to \{-1,1\}$. We are interested in Tree (T) and Tree Ensemble (TE) models, and in this section, we review the basic concepts of Ts and TEs that are needed in the following Sections. 

Decision trees yield a partition of the feature space $[0,1]^{|J|}$ by applying hierarchical disjunctive splittings.
A tree is characterized by a maximum depth $D$ so that nodes are organized into at most $D$ levels in $\D= \{0,1,\dots, D\}$. 
A decision tree is composed of \emph{branch nodes} in the set $\B$ and \emph{leaf nodes} in the set $\L$. A branch node $t\in\B$ applies a splitting rule on the samples in the node, while a leaf node $\ell\in\L$ acts as a collector of samples. In each leaf node, the same class is assigned to all the samples contained within it using a simple rule, such as the majority vote.

The most common decision trees are univariate employing axis-aligned splits \mnote{\cite{Breiman1984CART}} \lnote{where the feature and associated threshold are chosen so as to decrease child nodes impurity.}

In recent years, 
multivariate (or oblique) decision trees have been proposed, which may involve multiple features per split \lcomment{ \cite{Bennett1992,Bertsimas2017OptimalClassification}}. Each branching rule is defined by hyperplane $h_t(x)=a_t^Tx^i+b_t$, where $a_t\in\mathbb{R}^{|\J|}$ and $b_t\in\R$ (where the apex $^T$ denotes transposition of the vector).
These multivariate splits are much more flexible than univariate ones, at the expense of lower interpretability. However, the problem of determining them is much more complicated, and it can be tackled using Mixed Integer Linear Programming. Of course, univariate models are a subclass of multivariate ones, thus in general, the splitting rule at a node $t\in\B$ can be represented as 
\begin{equation}\label{eq: routing rules}
\text{if }h_t(x^i) \begin{cases}
    \leq 0 & x^i\text{ follows the left branch of }t\cr
 >0               & x^i\text{ follows the right branch of }t
\end{cases}.
\end{equation}

The hierarchical tree structure uses hyperplane splits that recursively partition the feature space into disjoint regions, each of which corresponds to a leaf node in the tree. The obtained tree is then used to classify unseen data.
A root-to-leaf unique path (decision path), which is a conjunction of predicates, leads to the prediction $\hat{y}_T\in\{-1,1\}$ made by the tree.




We define a \emph{Tree Ensemble} TE as a collection of
tree estimators $\E$ with associated weights $w_e$, with $e\in\E$.
The TE decision function $F_{{TE}}:\mathbb{R}^{|\J|}\to \{-1,1\}$  is obtained as the weighted majority vote of the decision function of its trees (ties are usually broken
in favor of the smaller index).

\section{VITE: a hierarchical VIsualization tool for TE}\label{sec:vite}
The underlying idea of our graphical representation has been inspired by the definition of variable importance given by Ishwaran \cite{Ishwaran2007} and the definition of Minimal Depth \cite{Ishwaran2012}.
The importance of a variable is related to the level at which it is used as variable splitting, assuming that splitting at the root node or nodes at shallower tree levels has a greater impact on variable importance than those used at deeper levels. 
Indeed features that are used overwhelmingly at the shallower depths are those that decrease impurity the most and therefore play a greater role in the classification of the samples.
We aim to generalize this principle and consider all the times a feature is used in the trees, and not just the first appearance. This gives a hierarchical perspective on the role played by features in the trees that make up the TE.
Indeed, a feature can show up multiple times at lower levels in the tree after its initial appearance, which determines its Minimal Depth. On the other hand, it may not reappear after its initial appearance. We aim to gain insight into the role of features along all the levels of the trees. 


Drawing inspiration from this concept, we consider for each feature the percentage frequency with which it is selected, at each level, in the forest's trees.
We define the $j-$th feature \emph{level frequency} for each $d\in\D$ as the ratio between the number of times the feature $j$ is used at level $d$ of all the trees in $\E$  and the total number of nodes 
that effectively apply a split at level $d$, thus not accounting for the possible leaves appearing at level $d$ (pruned nodes).
In this definition, we must account for the type of TE considered. Indeed, as remarked in \cite{Tan2020}, in Random Forest each tree is generated by an identical and independent process and contributes equally to the prediction. However, this is no longer true in boosted trees where individual trees are obtained by a boosting process which makes trees not independent and with a different contribution to the overall prediction. This aspect can be managed by weighting the contribution of each tree similarly to the approach proposed in \cite{Tan2020} for the Proximity measure. Furthermore, each decision tree in a TE can be allowed to use only a limited random subset of features as candidates for splitting at each node \cite{ho1998random}. Again, in this case, a weight can be used to compensate for the bias as proposed for the Minimal Depht in \cite{Ishwaran2012}. 
To formally introduce the level frequency of a feature, we assume, for the sake of simplicity, that nodes in each tree are numbered according to breadth-first indexing (increasing from left to right at each level $d$), starting from the root node, which is numbered zero, so that $t\in\T=\{0,1,2,\dots, 2^{D+1}-1\}$.  
Let us define the indicator function for each $j\in\J, t\in \T, e\in \E,$ as
$$ \ds{1}(j,t,e)=\begin{cases} 1 & \text{if feature $j$ is used at node $t$ of the tree $e$}\cr 0 & \text{otherwise}
\end{cases}
$$
and for each tree $e\in\E$, let  $\B^e(d)$ be the set of nodes at level $d$ of the tree $e$ which effectively apply a splitting rule. Thus, we \lnote{introduce} the following definition.


\begin{definition}[Level Frequency]\label{def:level_frequency}
    The frequency $f_{d,j}$ of a feature $j$ at a level $d\in\{0,\dots, D-1\}$ is 
\begin{equation}\label{level_frequency}
f_{j,d} =  \frac{1}{\displaystyle\sum_{e\in\E} |\B^e(d)|} \sum_{e \in \E} w^e \sum_{t\in\B^e(d)} \ds{1}(j,t,e),
\end{equation}
where $w^e$ $ \forall e\in\E$ are non-negative weights that take into account 
(i) the probability of a feature being sampled as a candidate for the splitting rule in a level of the tree $e\in\E$ and (ii) the contribution of each tree $e\in\E$ in predicting the outcome.
\end{definition}

According to this definition, we can define the level frequency matrix of dimension respectively $ |\J| \times D$ with elements $\{f_{j,d}\}_{i\in \J, d\in \D}$ \lcomment{, rescaled in percentage}.
Each column of the matrix sums up to 100\%.

We represent this matrix with a heatmap where the darker colors represent the more used features at the level $d$.

To obtain a deeper view of the features' use in the TE we can consider a more specific definition of the frequency for each node $t$ in the tree. 
We define the $j-$th feature \emph{node frequency} for each $t\in\T$ as the ratio between the number of times the feature $j$ is used at node  $t$ in all the trees in $\E$  and the total number of nodes $t$ that are splitting nodes in the ensemble. 
As in the definition \ref{def:level_frequency}, we can weigh the contribution
of each tree to account for a possible different role in the TE. 
Thus, we have the following definition.

\begin{definition}[Node Frequency]\label{def:node_frequency}
    The frequency $f^{\text{node}}_{t,j}$ of a feature $j\in\J$ at a node $t\in\T$ is 
\begin{equation}\label{node_frequency}
f^{\text{node}}_{t,j} =  \frac{1}{\displaystyle\sum_{e\in\E}\sum_{j\in\J}\ds{1}(j,t,e)} \sum_{e \in \E}w^e\ds{1}(j,t,e),
\end{equation}
where $\sum_{e\in\E}\sum_{j\in\J}\ds{1}(t,j)$ is equal to the number of times the node $t$ is a splitting node in the forest and $w^e$ $ \forall e\in\E$ are non-negative weights that take into account 
(i) the probability  of a feature being sampled as a candidate for the splitting rule in a node of the tree $e\in\E$ and (ii) the contribution of each tree $e\in\E$ in predicting the outcome.
\end{definition}

We use \lnote{the node frequency} definition to obtain a visualization of the TE in the form of a single tree of depth $D$ where for each node $t\in\T$ we report the heatmap of the features' node frequency matrix with elements 
$\{f^{\text{node}}_{t,j}\}_{i\in \J}$\lcomment{, rescaled in percentage}. Also in this case, for each $t\in\T$ we have $\sum_{j\in\J}f^{\text{node}}_{t,j}$ \gnote{$\times 100$ }$ =100\%$.
Of course, the leaf nodes which are at level $D$ are not represented since there are no splits. However, in case a terminal node appears in upper levels, from $0$ to $D-1$, it is reported with values $f^{\text{node}}_{t,j}=0$ for all $j\in\J$. \\
As additional information, we recover from the TE  the ranges of the threshold $b_j$ for each feature $j$ used in the trees of the TE at each node $t\in\T$. This is reported in the visualized tree, as an interval $[l_j,u_j]$ near the name of each feature $x_j$ which is used at node $t$. When $l_j=u_j$, a single value is reported. The interval gives us additional information: for each tree in the TE, whenever a features $j$ is used at node $t$ we trivially have that
$$
x^i_j\begin{cases}
 \le l_j & x^i \text{ always goes to the left branch}\cr
\in (l_j,u_j] & \text{uncertain (grey choice)}\cr
 > u_j & x^i\text{always goes to the right branch}
\end{cases}
$$
This can give additional insight that can be analyzed by experts in the field under study. Indeed it may allow identifying easily understandable decision paths for samples $x^i$ whose features $j$  fall into the external intervals of the $(-\infty,l_j]$ or $(u_j,+\infty)$ for all $j\in\J$.

\lnote{It is worth mentioning how the proposed frequency metrics and \texttt{VITE} representations compare to existing metrics for computing feature importance in TE such as the Mean Decrease Impurity (MDI - Gini importance) and Mean Decrease Accuracy (MDA - permutation importance), defined in \cite{breiman2001random}.  
 The features that most effectively reduce node impurity (and hence have higher MDI) are selected early and more frequently in the tree-splitting process, aligning with the underlying greedy algorithm for constructing decision trees in a forest. 
 Further, the features selected first and most frequently are decisive for evaluating MDA.  
Thus, the most important features, according to MDI and MDA, align closely with those identified through \texttt{VITE}, which gives similar insight into the overall features' ranking. 
However, traditional metrics, such as MDI and MDA, provide a ranking of features based on the average
role they play along the trees in the forest. Instead, \texttt{VITE} provides a visual representation that evaluates the importance of features
in a distributed manner along different levels of trees throughout the forest, providing a different hierarchical view. Moreover, our
visualization incorporates the threshold ranges that are applied in the splitting rules, offering insights into how specific features are
utilized at various depths within the forest. It is worth noting that the MDI metric can be adapted to evaluate the decremental
impurity contribution of each feature at various levels. Specifically, the vanilla MDI for a given feature is calculated as the total
sum of impurity reductions across all nodes in each tree where the feature is used, normalized by the total number of trees in the
ensemble. An interesting approach could be to evaluate the MDI for nodes at a certain level across all trees where the feature is
employed by the TE, utilizing this measure in the \texttt{VITE} representation instead of the frequency.
} 

\begin{figure}[!htb]
\centering
\includegraphics[width=1\textwidth]{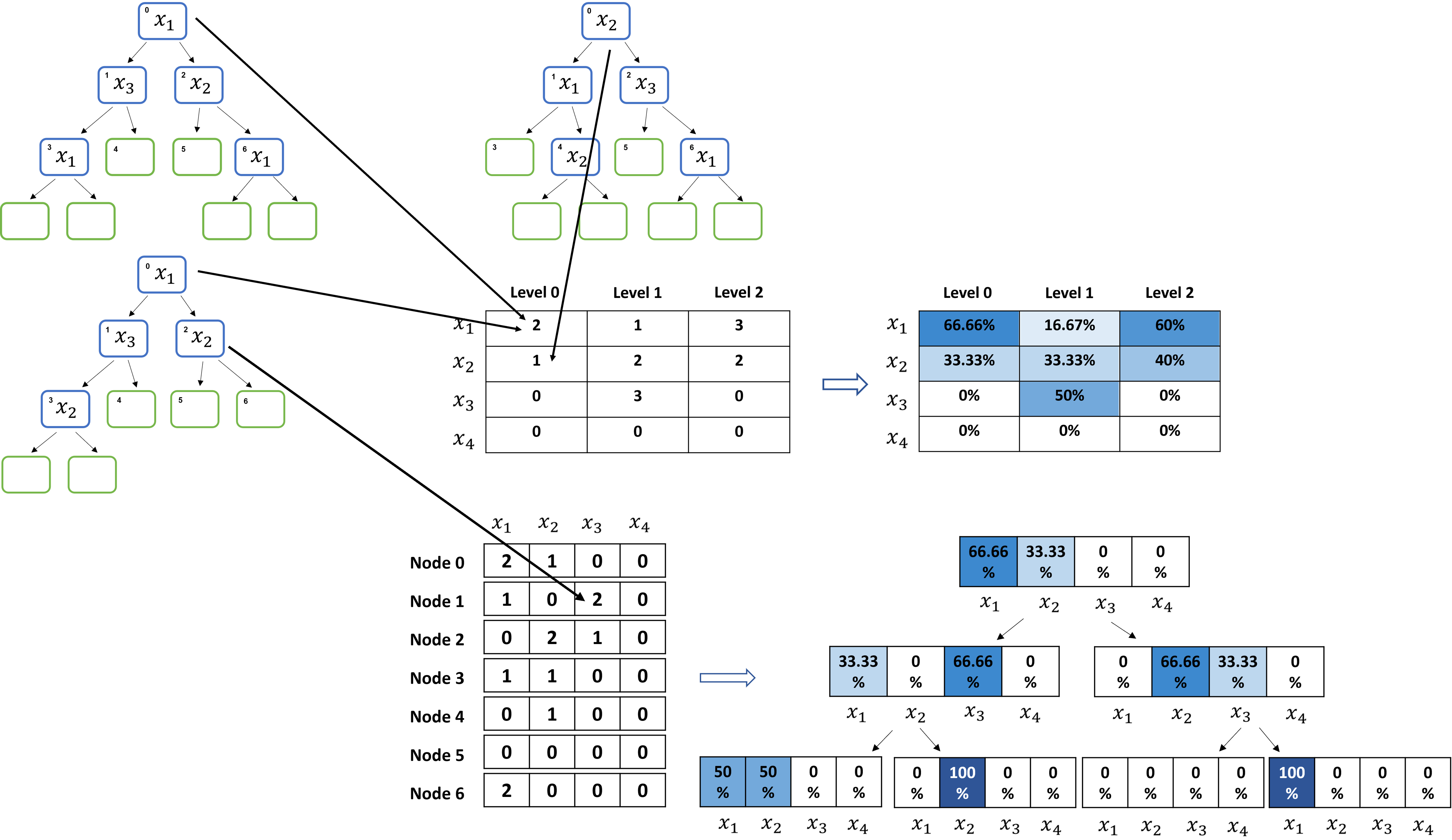}
\caption{\lcomment{Toy example: } construction of the features' usage heatmap at different depths in the TE}
\label{fig:visualization1}
\end{figure}

\gnote{We show on a toy example how \texttt{VITE} works} for the construction of the \gnote{heatmap} matrices and of the representative tree in Figure \ref{fig:visualization1}. 
Here, we assume to have a simple Random Forest composed of $E=3$ trees with depth $D=3$, and samples $x^i\in \R^4$, namely $\J = \{1,2,3,4\}$. 

\begin{figure}[!htb]
\centering
\includegraphics[width=.7\textwidth]{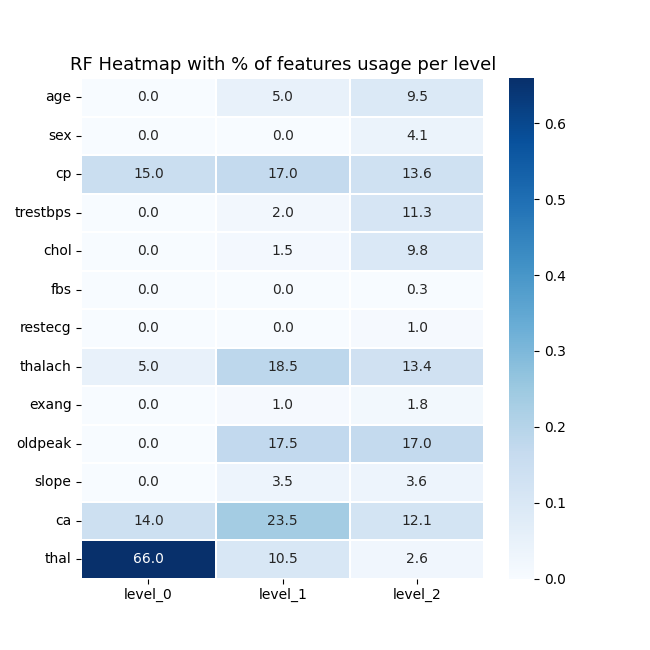}
\caption{\lcomment{Cleveland  example: features' level frequency at the three different tree levels of a RF}}
\label{fig:heatmap}
\end{figure}

\begin{sidewaysfigure}[!htb]
\centering
\includegraphics[width=1.12\textwidth]{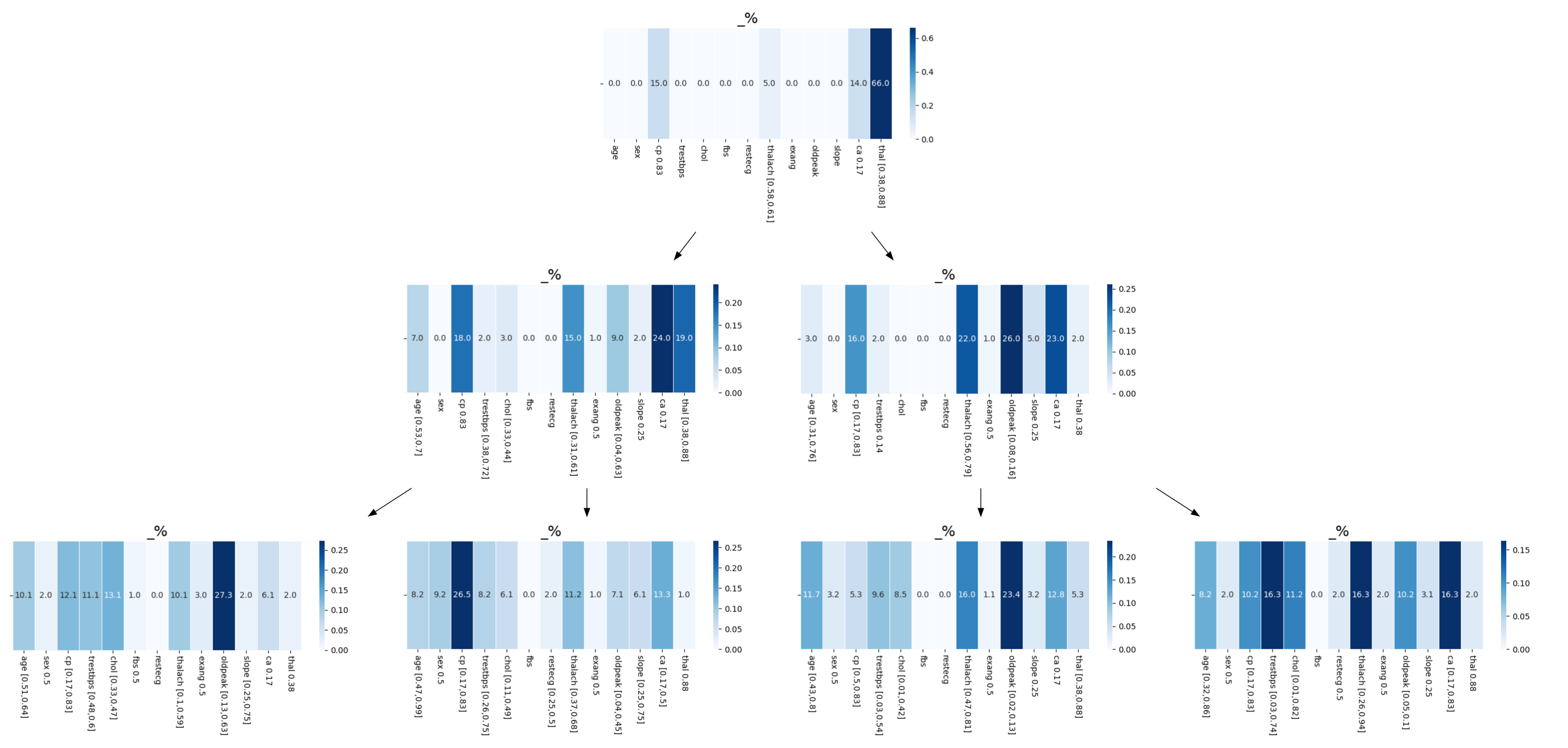}
\caption{\lcomment{Cleveland  example: features' node frequency at the nodes of the trees in the RF}}
\label{fig:tree_heatmap}
\end{sidewaysfigure}

\begin{figure*}[ht!]
  \centering
  \begin{subfigure}[b]{0.45\textwidth}
    \includegraphics[width=\textwidth, height=4cm]{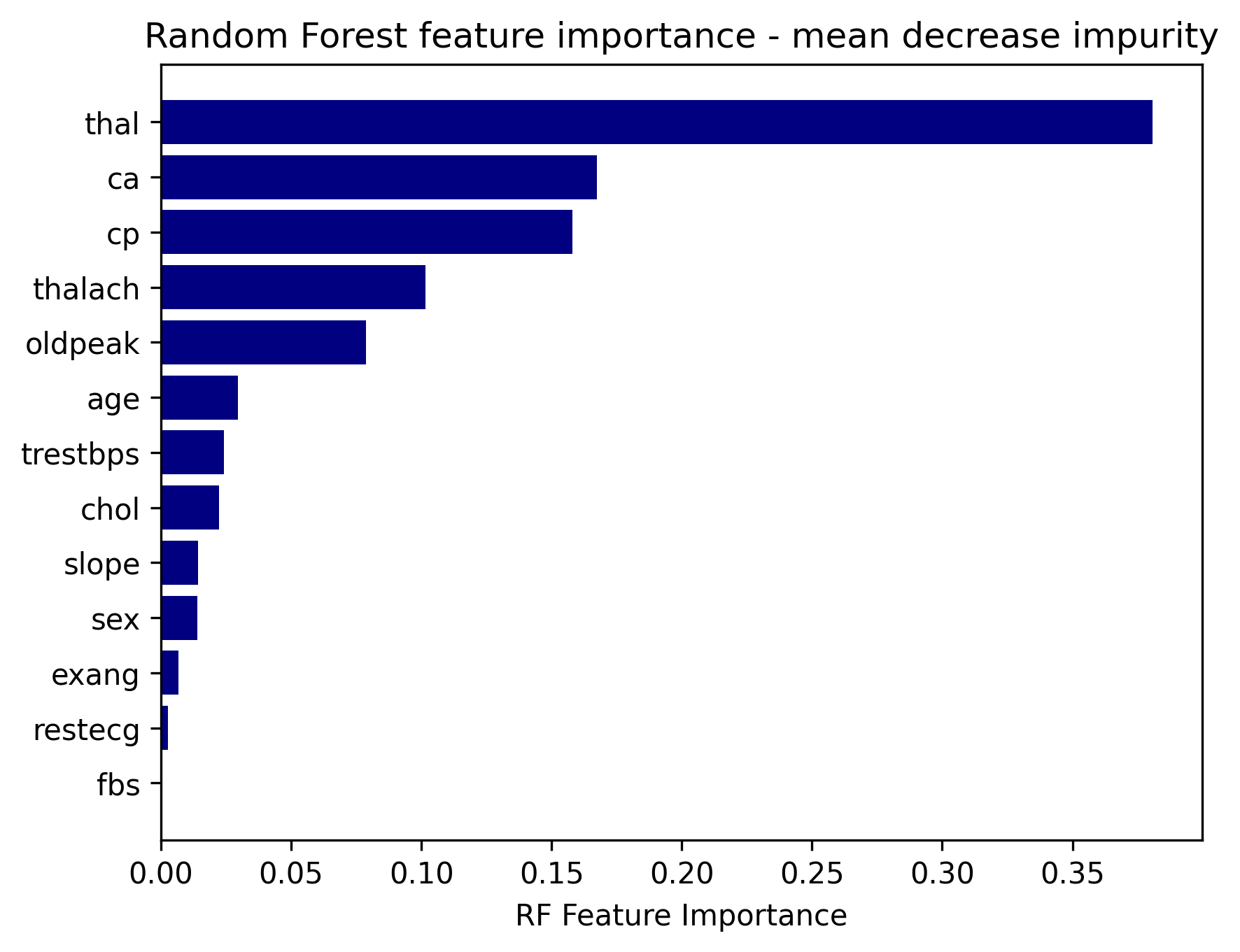}
    \caption{}
    \label{fig:MDI}
  \end{subfigure}
  \hfill
\begin{subfigure}[b]{0.45\textwidth}
    \includegraphics[width=\textwidth, height=4cm]{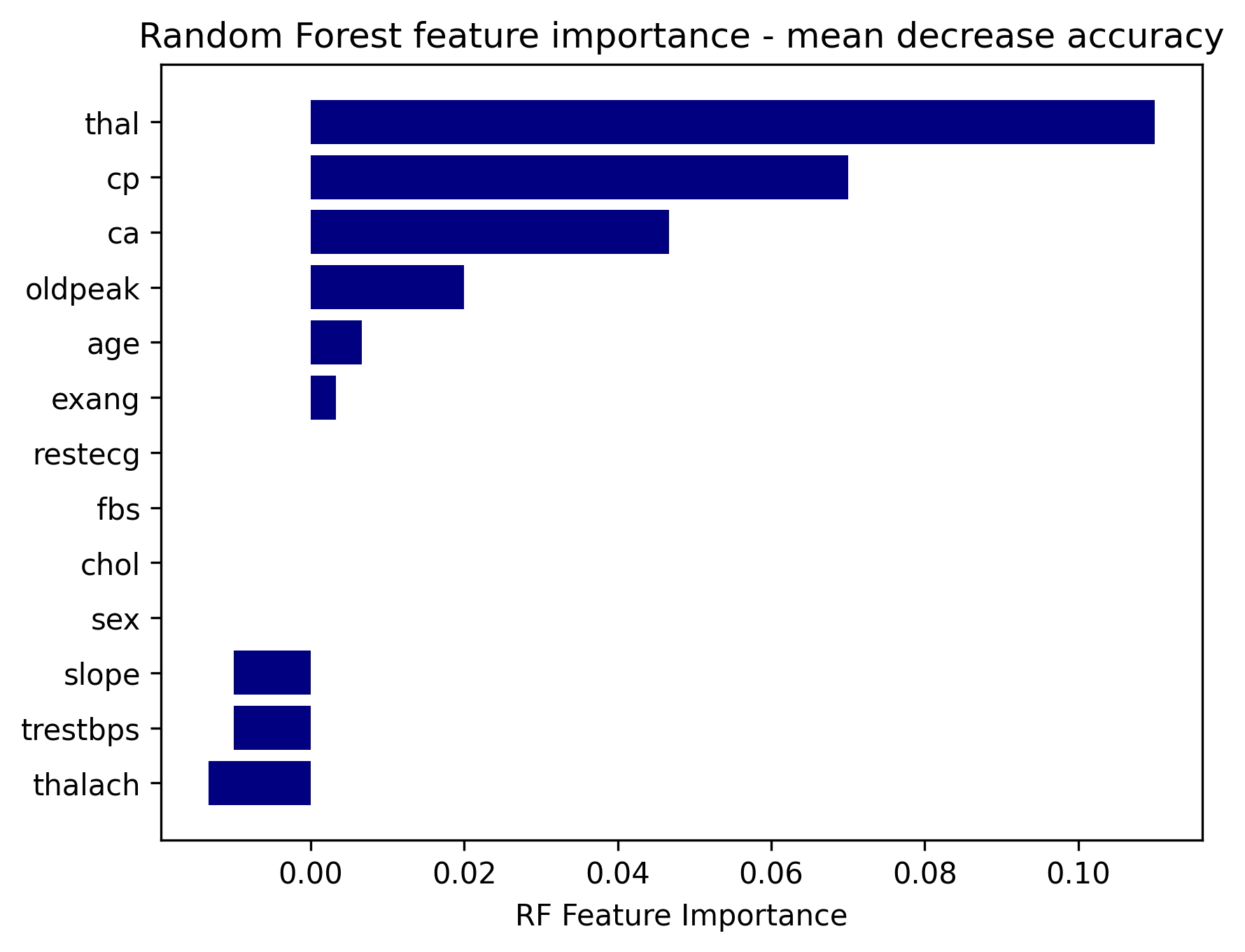}
    \caption{}
    \label{fig:MDA}
  \end{subfigure}
\caption{\lcomment{Cleveland example: Figure (a) represents the features' importance by mean decrease impurity (MDI). Figure (b) represents the features' importance by mean decrease accuracy (MDA).}}
  \label{fig:feature_importance_cleveland}
\end{figure*}

\lcomment{In this simple example, features $x_1$ and $x_2$ are both used at the root node and stand out as the most frequently used in the TE, also at shallower depths. Feature $x_3$ appears exclusively at the $d = 1$ level, while feature $x_4$ remains unused in the entire forest. From such a representation it would seem reasonable to infer that features $x_1$ and $x_2$ are critical for impurity reduction when most samples are accessible to a node, and consequently for predictions. Feature $x_3$ is less important and $x_4$ is not effective as a predictor. The intuitive nature of these visualizations provides simple insights into feature usage, making the model understandable even for non-expert users.\\}
\lnote{We comment on the \texttt{VITE} representation, and the insight gained when applied for  the visualization of a Random Forest  with $E=100$ tree estimators of maximum depth $D=3$ trained on Cleveland dataset \cite{UCI2019}, which refers to heart disease classification and has  $x\in \R^{13}$. In this case, features have been normalized in $[0,1]$, and we disable random sampling of the features in the construction of the RF. We use weights $w^e=1$ in the definitions of frequencies \eqref{def:level_frequency} \eqref{def:node_frequency}. In Figure \ref{fig:heatmap}, we report the heatmap representing the level frequency matrix; in Figure \ref{fig:tree_heatmap}, we report the representative tree of depth $D=3$ where at each node $t$, the node frequency of each feature $j$ is reported together with the interval $[l_j,u_j]$.
In Figure \ref{fig:feature_importance_cleveland}, we report the aggregated metrics MDI and MDA to compare with.
 From \texttt{VITE} visualization, only four features ("thal", "cp", "ca" and "thalach") are used at all the levels of the trees with different percentages, being  "thal" the one with the prominent role (66\%) at the root node. 
 Comparing with rankings of Figure 
 \ref{fig:feature_importance_cleveland} derived from the MDI   and MDA, it is evident, as expected, that the most important features "thal", "cp", "ca",
 align closely with those identified by \texttt{VITE}.
 In contrast to the aggregated perspective of these traditional metrics, \texttt{VITE} offers a distributed representation and a view of the threshold ranges used at different levels. For example,  "thal" is the most prominent feature used in the root node, confirming as known from the medical literature, that the Thallium stress test result (thal) is a strong predictor of potential heart disease.
Further, splitting values for "thal" are always in the same range $[0.38, 0.88]$ at each level of all the trees, giving insight into the values used to discriminate among patients in the nodes.
In a hierarchic view, at the next level, "oldpeak" and "age" come into play at a significant percentage (greater than 5\%).  
On the other hand, some features, e.g. "sex" which is binary, appear in the TE trees only at the very last branching nodes, thus showing a role only as a final discriminator.}

The source code and the data to use the \texttt{VITE} tool are available at \url{https://github.com/gditeodoro/VITE}.

\clearpage

\section{\texttt{MIRET}: a Multivariate Interpretable REbuilt optimal Tree}
\label{sec:model}
\subsection{Introduction}
In this section,  we proposed a Mixed Integer Linear Programming (MILP) formulation for learning a single optimal multivariate tree as a surrogate model of the ensemble classifier.

In more detail,  let ${TE}$ be a \textit{target Tree Ensemble model}, with ${|\E|}$ tree estimators of maximal depth $D$ trained on the given dataset and let $F_{{TE}}:\mathbb{R}^{|\J|}\to \{-1,1\}$  be the decision function of  ${TE}$.

The aim is to provide a shallow and interpretable surrogate multivariate tree model $T$, of the same maximal depth $D$ of ${TE}$, with 
decision function
$F_T: \mathbb{R}^{|\J|}\to \{-1,1\}$,
which is able to reproduce as much as possible the ${TE}$ predictions on the training set.
Further, we aim to enforce in the splitting hyperplanes the use of the most informative features.

As a first step, we formalize the definition of \emph{fidelity} we are using. 
\begin{definition}[Training data fidelity]
Given a tree ensemble ${TE}$, we say that a multivariate decision tree $T$
of depth $D$ is faithful to training data w.r.t.
 ${TE}$,
when the prediction 
$F_{{TE}}(x^i)=F_T(x^i)$ for all $i\in\I$.
\end{definition}
This definition is different from the one proposed in \cite{vidal2020born}, where instead, authors construct a  \emph{whole faithful}  representer tree which aims to reproduce exactly the decision function on the whole space of the features, namely  $ F_{{TE}}(x)=F_{{T}}(x) $, {for all $x\in\R^{|\J|}$}, and not only on the training samples $x^i, {i\in\I}$, or from the one given in \cite{breiman2001random} where manufactured data are used in the definition of fidelity. 
In those approaches, the complexity (depth and/or number of splits) is minimized, exactly or heuristically, but in principle, it can grow exponentially.
Instead, in our approach, we aim to find a partition of the space which produces the same classification on the available samples keeping the depth of the representer tree fixed to $D$, the same as the trees in the $TE$. Since depth $D$ is fixed, the existence of such a faithful $T$ on the training data is not straightforward and, in general, cannot be guaranteed.
Thus, we chose to maximize the fidelity and {to include a penalization term to favor sparsity in the splits}. According to the taxonomy reported in \cite{burkart2021survey}, we refer to the following definition of 
 \emph{global surrogate fitting model}.
Giving a loss function $S:\{-1,1\}^{|\J|}\times \{-1,1\}^{|\J|} \to \R$ which measures the error in binary classification {and a penalty function $V(T)$, which measures the number of features used by the tree $T$}, we aim to find a tree $T^*$ in a restricted set of decision trees  ${\cal F}_T$ such that
$$T^*=\arg\min_{T\in {\cal F}_T} \left(\sum_{i\in \I} S(F_T(x^i),F_{{TE}}(x^i))+ V(T)\right)  $$ 
In our approach, trees in the class
${\cal F}_T$ are characterized by
\begin{itemize}
    \item fixed depth $D$;
    \item branching at node $t$ according to rule \eqref{eq: routing rules} with a multivariate splitting \begin{equation}\label{eq:hyperplane}
h_t(x)=a_t^Tx+b_t,   
\end{equation}
 with $a_t\in\R^n$ and $b_t\in\R$;
 \item constraints derived from $TE$ information.
\end{itemize}

{The training phase outputs an optimal classification tree $T^*$ defined by coefficients $a^*_t$ and $b^*_t$ for each branching node $t\in\Tb$.}

\smallskip
\paragraph{Notation}
In order to present our formulation, we introduce the main concepts and notation concerning optimal trees, as first described in \cite{Bertsimas2017OptimalClassification}.
Since the tree has fixed depth $D$,
 levels can be numbered as $\D=\{0,\dots,D\}$ being $0$ the root node and $D$ the terminal level of the leaf nodes.
As usual, nodes are divided into \emph{branch} and \emph{leaf} nodes. 
A branch node $t$ applies the multivariate splitting rule \eqref{eq:hyperplane} on the subset of samples $\I_t\subseteq \I$ assigned to it, partitioning them among the left or right branch and hence among the two child nodes according to \eqref{eq: routing rules}.
Following \cite{MarginOptimalClassificationTrees}, branching nodes $\B$ at level $d<D$ always apply a split which can be an effective or a dummy one. Indeed, if a node $t$ at a level $d < D$ does not need to partition the samples further, we define two dummy children anyway such that only one of the two contains all the samples $\I_t$. 
We assume that nodes in the tree are numbered according to breadth-first indexing, starting from the root node, which is numbered zero and increasing from left to right at each level $d$, so that the tree has an overall number of nodes $\displaystyle\sum_{k=0}^D 2^k-1= 2^{D+1}-2$. 
      
\begin{figure}[ht]
\centering
\hspace{1.2cm}
\scalebox{.45}{{\tikzset{
  branchnode/.style = {line width=1.8pt, scale =1.2, very thick,          rectangle,
        rounded corners, draw=dblue, text width=2em, text centered, minimum height=1cm, minimum width = 1cm, font=\Large},
  1dummy/.style       = {font=\Large},
  rect/.style = {draw = none, font=\LARGE, text width = 5 cm},
    rect2/.style = {draw = none, font=\LARGE}
}

\tikzset{
    ncbar angle/.initial=90,
    ncbar/.style={
        to path=(\tikztostart)
        -- ($(\tikztostart)!#1!\pgfkeysvalueof{/tikz/ncbar angle}:(\tikztotarget)$)
        -- ($(\tikztotarget)!($(\tikztostart)!#1!\pgfkeysvalueof{/tikz/ncbar angle}:(\tikztotarget)$)!\pgfkeysvalueof{/tikz/ncbar angle}:(\tikztostart)$)
        -- (\tikztotarget)
    },
    ncbar/.default=0.3cm,
}

\tikzset{square left brace/.style={ncbar=0.3cm}}
\tikzset{square right brace/.style={ncbar=-0.3cm}}

\begin{tikzpicture}
  [
    grow                    = down,
    level 1/.style          = {sibling distance=11.5cm},
    level 2/.style          = {sibling distance=5.7cm},
    level 3/.style          = {sibling distance=3.2cm},
    level distance          = 3.6cm,
    edge from parent/.style = {thick, draw, edge from parent path={(\tikzparentnode) -- (\tikzchildnode)}},
  ]
\node [branchnode] (0) {0}           
    child {
    node[branchnode] (1) {1}
    child { node [branchnode] (3) {3} 
    child { node [branchnode, draw = verdone] (7) {7}
         edge from parent node[1dummy, pos=0.55, left=0.4cm] (e) {$\leq$} }
    child { node [branchnode, draw = verdone] (8) {8}
      edge from parent node[1dummy, pos=0.55, right=0.4cm] (k) {$>$}}
    edge from parent node[1dummy, pos=0.4, left=0.55cm] (c) {${\dblue{a}}^T_1 x^i + \dblue{b}_1 \leq 0$}}
    child { node [branchnode] (4) {4} 
    child { node [branchnode, draw = verdone] (9) {9}
    edge from parent node[1dummy, pos=0.55, left=0.4cm] (h) {$\leq$} }
    child { node [branchnode, draw = verdone] (10) {10}
      edge from parent node[1dummy, pos=0.55, right=0.4cm] (k) {$>$}}
    edge from parent node[1dummy, pos=0.4, right=0.55cm] {$>$} }
  edge from parent node [1dummy, pos=0.35, left=0.85cm] (a) {${\dblue{a}}^T_0 x^i + \dblue{b}_0 \leq 0$}
}
    child {
     node[branchnode] (2) {2}
    child { node  [branchnode] (5) {5} 
    child { node [branchnode, draw = verdone] (11) {11}
        edge from parent node[1dummy, pos=0.55, left=0.4cm] (j) {$\leq$} }
    child { node [branchnode, draw = verdone] (12) {12}
        edge from parent node[1dummy, pos=0.55, right=0.4cm] (k) {$>$} }
     edge from parent 
     node [1dummy, pos=0.4, left=0.55cm] {$\leq$}}
    child { node [branchnode] (6) {6} 
    child { node [branchnode, draw = verdone] (13) {13}
      edge from parent node[1dummy, pos=0.55, left=0.4cm] (l) {$\leq$}}
    child { node [branchnode, draw = verdone] (14) {14}
        edge from parent node [1dummy, pos=0.55, right=0.4cm] (f) {$>$}}
    edge from parent 
     node[1dummy, pos=0.4, right=0.55cm, draw=none] (d) {$\dblue{a}^T_2 x^i + \dblue{b}_2 > 0$}}
    edge from parent node [1dummy, pos=0.35, right=0.85cm] (b) {${\dblue{a}}^T_0 x^i + \dblue{b}_0 > 0$}
};





\node[rect,  right = 6.8 cm of b] (Tbf) {\dblue{$\Tbf=\B(0)\cup\B(1)$}};
\node[rect,  below = 4.7 cm of Tbf] (Tbl) {\dblue{$\Tbl=\B(2)$}}; 
\node[rect, below = 8.7 cm of Tbf] {\verdone{$\Tl$}}; 


\draw [black, thick] (12,0.8) to [square left brace] (12,-4.2);
\draw [black, thick] (12,-5.8) to [square left brace] (12,-8.3);
\draw [black, thick] (12,-9.8) to [square left brace] (12,-12.3);


\draw [black, thick, rotate = -90] (12.5,-6.0) to [square left brace] (12.5,-11.4)  ;

\draw [black, thick, rotate = -90] (12.5,-0.1) to [square left brace] (12.5,-2.5);
\draw [black, thick, rotate = -90] (12.5,-3.3) to [square left brace] (12.5,-5.7);

\draw [black, thick, rotate = 90] (-12.5,-11.4) to [square left brace] (-12.5,-0.1);


\node[1dummy, below = 0.03cm of 7]  {$c_7=-1$};
\node[1dummy, below = 0.03cm of 8]  {$c_8=1$};
\node[1dummy, below = 0.03cm of 9]  {$c_9=-1$};
\node[1dummy, below = 0.03cm of 10]  {$c_{10}=1$};
\node[1dummy, below = 0.03cm of 11]  {$c_{11}=-1$};
\node[1dummy, below = 0.03cm of 12]  {$c_{12}=1$};
\node[1dummy, below = 0.03cm of 13]  {$c_{13}=-1$};
\node[1dummy, below = 0.03cm of 14]  {$c_{14}=1$};

\node[rect2,  below = 5.3cm of 3] (a) {\verdone{$\S_L(1)$}};
\node[rect2, right = 12.6 cm of a]  {\verdone{$\S_R(0)$}};

\node[rect2, right = 2.3 cm of a]  {\verdone{$\S_L(4)$}};
\node[rect2, right = 5.3 cm of a]  {\verdone{$\S_R(4)$}};

\end{tikzpicture}}}
\caption{Representation of structure and notation for a tree with depth $D=3$.}
\label{fig: tree_graph}
\end{figure}
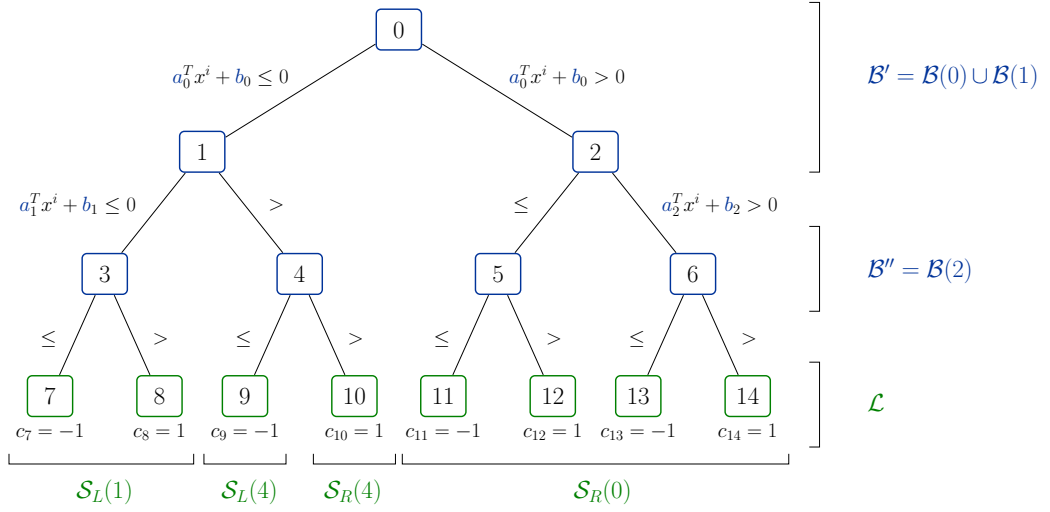

We  make use of the following notation, which is pictured in Figure \ref{fig: tree_graph}:
    \begin{itemize}
    \item $\Tb$: the set of \textit{branch nodes}
     where an oblique splitting rule is applied; $\Tb$ are numbered $\{0,\dots, 2^D-2\}$. 
    \item $\Tl$: the set of \textit{leaf nodes} where a class is assigned to a sample; 
    $\Tl$ are  numbered $ \{2^D - 1  , \dots , 2^{D+1} - 2\}$.
    \item $\B(d)$ the set of  nodes at level $d$, numbered as  $\{0,\dots, 2^D-2\}$. ;
    \item $\displaystyle \Tbf=\bigcup_{d=0}^{D-2}\B(d)$: the set of branch nodes not adjacent to the leaves;
        \item $\Tbl =\B(D-1)$: the set of branch nodes adjacent to the leaves; 
\item $\S(t)$ the set of leaf nodes in the subtree rooted at node $t\in\B$ 
\item  \lnote{$ \S_L(t)$ and $ \S_R(t)$ the set of left and right sub-leaves respectively,} i.e. the set of leaf nodes following the left  and right branch of the subtree rooted at node $t\in\B$;
    \item the class assignment $c_\ell$ for each $\ell \in \L$. Being in a binary classification setting, we pre-assigned a class label to each leaf node, labeling as $-1$ the odd leaves and as $+1$ the even ones. 
    \end{itemize}

\paragraph{Variables}
The straightforward variables are the coefficient of the hyperplane $h_t$ for each $t\in\B$ which are continuous and w.l.o.g. can be assumed to be normalized so that
$$a_t=\{a_{t,j}\}_{j\in \J}\in[-1,1]^{|J|}, \  b_t\in[-1,1] \quad \forall \ t\in\B.$$




As done in \lnote{traditional} OCTs \cite{Bertsimas2017OptimalClassification}, 
we also need the binary variables $ z_{i,\ell} $ for all samples $i\in \I$ and leaves $\ell\in \Tl$, defined as follows:

$$
\begin{array}{llll}
    z_{i,\ell} =
    \begin{cases} 
      1 & \text{if sample $i$ is assigned to leaf $\ell\in \Tl$} \\
      0 & \text{otherwise} 
   \end{cases}.
\end{array}
$$


Since we are interested in controlling the use of features along the tree, we  introduce other binary variables $s_{t,j}\in\{0,1\}$, for each $t\in\Tb$ and $j\in\J$, as follows:

$$
\begin{array}{llll}
    s_{t,j} =
    \begin{cases} 
      1 & \text{if feature $j$ is \emph{used} at node $t$ ($a_{t,j}\neq 0$)} \\
      0 & \text{otherwise}
   \end{cases}.
\end{array}
$$

Table \ref{tab:variables_overview} presents an overview of the variables used in the model.

\begin{table}[ht!]
\renewcommand\arraystretch{1.4}
\footnotesize
\centering
\begin{tabular}{lll}
\textbf{Variable}                & & \textbf{Description}                                             \\ \hline\hline
$a_{t,j} \in [-1,1]$    & & hyperplane coefficient at node $t$ of feature $j$       \\
$b_{t} \in[-1,1]$     & & hyperplane intercept at node $t$               \\
$s_{t,j} \in \{0,1\}$   &  &1 if feature $j$ is selected at node $t$, 0 otherwise  
\\
$z_{i,\ell} \in\{0,1\}$ &  &1 if sample $i$ is assigned to leaf $\ell$, 0 otherwise \\
\end{tabular}\caption{Overview of decision variables in {\texttt{MIRET}}}
    \label{tab:variables_overview}
\end{table}

\paragraph{Tree structure based Constraints}
{We need to state constraints to recover the tree structure  as in MILP-based OCT formulation \lnote{for multivariate trees} (e.g. \cite{Bertsimas2017OptimalClassification}, \cite{MarginOptimalClassificationTrees}, \mnote{\cite{boutilier2022shattering}}).
}

The first set of constraints forces each sample $x^i$, $i\in\I$ to be assigned to one and only one leaf node $\ell\in\L$.
\textit{Assignment constraints}  are stated as follows:

\begin{align}
& \sum_{\ell\in\Tl} z_{i,\ell} = 1   && \forall\  i\in\I. \label{cons: fin assignment}
\end{align}



We further need to model disjunctive conditions on samples that model the routing rules defined in  \eqref{eq: routing rules} and ensure that the hyperplane splits are designed accordingly to the assignment $z$ of samples to the leaf nodes. 
The  \textit{routing constraints} on the sample $x^i$ are defined at each node $t\in\B$ but
must apply only to that $ i\in\I_t$, which is determined during the optimization process itself.

As in \cite{MarginOptimalClassificationTrees}, we use the observation that whenever $i\in \I_t$,  then $x^i$ must end up in one of the leaves of the subtree rooted at  $t$, and more specifically, either in the subset of the left or of the right subtree rooted
at node $t$. This condition is expressed as
$$\text{either  }\sum_{\ell\in\S_L(t)} z_{i,\ell}=1 \qquad \text{or}\quad \sum_{\ell\in\S_R(t)} z_{i,\ell}=1 .$$ 
Thus, we can write the \emph{routing constraints} as 
 
\begin{align}
  & \displaystyle a_t^T x^i + b_t \leq M_{\H_L} \Big (1 - \sum_{\ell\in\S_L(t)} z_{i,\ell} \Big) && \forall i\in\I ,\   \forall t\in \Tb \label{cons: routing left},
\\
& \displaystyle a_t^T x^i + b_t - \varepsilon  \geq -M_{\H_R} \Big (1 - \sum_{\ell\in\S_R(t)} z_{i,\ell}\Big) && \forall i\in\I,  \ \forall t\in \Tb
 \label{cons: routing right}
\end{align}
where $ \varepsilon >0 $ is a sufficiently small positive value to model strict inequalities in \eqref{eq: routing rules}. \lnote{In our implementation, we fixed $\varepsilon = 0.001$ (see \cite{Bertsimas2017OptimalClassification}).}

Indeed, if a sample $i\not\in \I_t$, we have that $\sum_{\ell \in \S_L(t)} z_{i,\ell} = 
\sum_{\ell \in \S_R(t)} z_{i,\ell}=0$, and thus both the constraints do not impose any routing conditions on sample $i$.
Instead, if sample $i\in\I_t$ and $\sum_{\ell \in \S_L(t)} z_{i,\ell} = 1$, thus, only the constraint \eqref{cons: routing left} is activated for sample $i$ at node $t$, while \eqref{cons: routing right} is deactivated. Analogous considerations can be done if a sample $i\in\I_t$  follows the right branch.

The Big-M values can be easily obtained because $x^i\in [0,1]^{|\J|}$, $ a_t\in[-1,1]^{|\J|}$ and $b_t\in[-1,1]$  for all $t\in\Tb$. Hence we can set the Big-M values to 
{$$M_{\H_L}=|\J|+1 \qquad M_{\H_R}=|\J|+1+\varepsilon.$$}
On the other hand, the choice of the $\varepsilon$ parameter is critical because small values may lead to numerical issues, and large ones may cut feasible solutions.

{Additionally, we need to include constraints among $s_{t,j}$ and $a_{t,j}$ to force the conditions $s_{t,j}=0 \ \Longleftrightarrow \ a_{t,j}=0$. This is easily modelled by \emph{sparsity constraints} as
\begin{equation}\label{cons: sparsity}
-s_{t,j} \leq a_{t,j} \leq s_{t,j}.
\end{equation}







\subsection{Incorporating TE-driven information}
We embed information derived from the $TE$ both in the constraints and in the objective function of the model \texttt{MIRET}.
In particular, we use
\begin{itemize}
    \item the proximity measure $m_{ik}$ among samples $x^i,x^k$; 
    \item the frequency $f_{d,j}$ of each feature $j$ at each level $d$;
    \item the probability $p^i$ of the $TE$ class prediction $\widehat y^i$,
    
\end{itemize}
and we report here the definitions needed.

The \textit{proximity measure}  $m_{ik}$ of a pair of samples $x^i,x^k$ in the $TE$ is
the number of trees $e\in\E$ in the $TE$ model in which
the two samples end up in the same leaf.
Let us define the indicator function for each pair $i,k\in \I$ and each tree $e\in\E$ as
$$ \ds{1}(i,k,e)=\begin{cases} 1 & \text{if samples $x^i, x^k$ end in the same leaf $\ell$ of the tree $e$}\cr 0 & \text{otherwise}
\end{cases}
$$
 Thus, we have the following definition \cite{breimancutler,Tan2020}:
  
\begin{equation}\label{eq:proximity}
m_{i,k} = \frac{1}{|\E|} \displaystyle\sum_{e \in \E} w^e \sum_{e\in\E} \ds{1}(i,k,e),
\end{equation}
where $w^e$ $ \forall e\in\E$ are non-negative weights that play a similar role as in the definitions \ref{level_frequency} and \ref{node_frequency}.

The proximity measure is considered a measure of the distance in data space. Two samples being in the same leaf in all the trees indicates that the two samples always belong to the same data space partition. We use it in the constraints. 
\par\medskip

The \emph{level frequency} $f_{d,j}$ of each feature $j$ at each level $d$ across the trees in $TE$ is calculated according to \eqref{level_frequency}. 

In our optimization model, we grow a full tree with depth $D$, and all the branch nodes at level $d$ are potential splitting nodes. In total, nodes are
$|\E| \cdot 2^d$, where $2^d$ is the number of nodes at level $d$. Since we aim to measure the spread of features among nodes in the trees, we use in  \eqref{level_frequency} the denominator equal to $|\E| \cdot 2^d$, and we have for all $d\in\D$, $j\in\J$:

\begin{align}
&f_{d,j} = 
\displaystyle\frac{1}{|\E| 2^d} \sum_{e \in \E} w^e \sum_{t\in\B^e(d)} \ds{1}(j,t,e). \label{eq: miret_frequency}
\end{align}
We use it both in the constraints and in the objective function.

\par\medskip
The \textit{class probability}  $p^i$  of a sample $x^i$ is the highest estimated probability of the prediction given by $TE$. It is calculated as 
\begin{align}
p^i := \ds \max_{c\in\{-1,1\}} p^i_c, \nonumber
\end{align}
where
\begin{align}
p^i_c := \frac{1}{E}\sum_{e=1}^E \frac{N_{\ell_e,c}}{N_{\ell_e}}, \nonumber
\end{align}
with $N_{\ell_e}$ the number of samples assigned to leaf $\ell$ in the tree  $e$ of the $TE$ model and as $N_{\ell_e,c}$ the number of samples with class label $c$ assigned to the same $\ell$ in $e$.

\lnote{The class probability metric derives from a procedure of soft voting of forest trees and has its roots in  \cite{breiman2001random}. 
 In this paper,  Breiman introduces the concept of 'votes' of individual forest trees, which are essentially used to determine the final label of a sample with a majority voting system.  
In majority voting, the predicted class label of a sample is the class label that represents the majority (mode) of the class labels predicted by each individual classifier.
This concept is then naturally extended to obtaining class probabilities with a soft voting of the forest's trees, as expressed in \cite{Tan2020}.
The predicted class probability of an input sample is the average of the predicted class probabilities of the trees in the forest. The class probability of a single tree is the fraction of samples of the same class in a leaf.
This is how sample class probability is implemented in \texttt{scikit-learn} library for Python \cite{scikit-learn}.}

We use it in the objective function.

\paragraph{TE-driven Constraints}

We use the \textit{proximity measure} to constrain the data space partitioning of the representative tree to be similar to the $TE$.
To this aim, we define   the set $\U_{TE}$ of pairs of samples with proximity measure higher than a given \textit{proximity threshold} $\overline{m}_{TE}\le 1 $ which is a hyperparameter of our model:
\begin{align}
&\U_{TE}{(\overline{m}_{TE})}:= \{(i,k)\in\I\times\I: \ i<k\  \wedge \  m_{i,k} {\ge} \overline{m}_{TE} \},\label{eq:proxiset}
\end{align}
and we include the \emph{proximity constraints}:
\begin{align}
&z_{i,\ell} = z_{k,\ell}  &&\forall \ell \in\Tl, \quad \forall (i,k)\in \U{(\overline{m}_{TE})} \label{cons: prox}
\end{align}
which ensures that pairs of samples with a proximity measure over the threshold are assigned to the same leaf in the model.


The \emph{level frequency} $f_{d,j}$ is used both to define hard constraints and to drive the selection of the features according to their   frequency.
Let us define the index sets
\begin{align} 
&\J_\gamma(d) := \{j\in\J: f_{d,j} > \gamma_d\}, && \forall d \in \D, \label{def:gamma_d}
\end{align}
where $\gamma_d$ for  $d\in\D$ are given \textit{frequency thresholds} and they are treated as hyperparameters of our model. 


The \emph{frequency constraints} force the model to select at level $d$ only features in the subset $\J_\gamma(d)$ by imposing the following conditions:
\begin{align}
& a_{t,j} = 0 && \forall j \in \J\setminus \J_\gamma(d)  && \forall t \in \Tb(d), && \forall d\in\D, \label{cons: freq}\\
& s_{t,j} = 0, && \forall j \in \J\setminus \J_\gamma(d)  && \forall t \in \Tb(d), && \forall d\in\D \label{cons: freq2}
\end{align}
This way, we induce local sparsity of the splitting hyperplanes in a hard way.\\

\smallskip


\subsection{Objective function}

The objective function is obtained as a combination of the two terms which measure fidelity and sparsity. 
In order to measure fidelity,  we need to define a loss function
 $S:\{-1,1\}^{|\J|}\times \{-1,1\}^{|\J|} \to \R$ to measure the error $S(F_T(x^i),F_{TE}(x^i))$.
The values $F_{TE}(x^i)=\widehat y^i$ are the class label predicted by $TE$ for sample $x^i$, based e.g. on majority voting.

The value $F_T(x^i)$ can be easily obtained as $$F_T(x^i)=\sum_{\ell\in\Tl} c_{\ell}z_{i,\ell}$$
where $c_\ell$ is the vector of classes assigned to leaves.
Thus $$\frac 1 2 \widehat{y}^{i}\Big(\widehat{y}^{i} - F_T(x^i)\Big)=\begin{cases}
    1 & \text{if $F_{TE}(x^i)=F_T(x^i)$} \cr
    0 & \text{otherwise}
\end{cases}.$$

We use class probability $p^i$ of the sample $x^i$ to weight misclassification with the aim of penalizing more the error on samples predicted with high probability by $TE$; we define the following loss function
$$S(F_{TE}(x^i), F_T(x^i) )=\frac 1 2 p^{i}F_{TE}(x^i)\Big(F_{TE}(x^i) - {F_T(x^i)}\Big)
=\frac 1 2 p^{i} \widehat{y}^{i}\Big(\widehat{y}^{i} - \sum_{\ell\in\Tl} c_{\ell}z_{i,\ell}\Big).
$$






We also aim to promote the use of sparser hyperplanes, enhancing the interpretability of the tree model.
Thus we penalize the selection of features  according to the level frequencies with which they are used at each level $d$ in $TE$. 
In particular, for each feature $j\in\J$ we count how many times it is used in the model at level $d$ as $\sum_{t\in\Tb(d)} s_{t,j}$. 
At each level $d$, the use of feature $j$ is weighted by the reciprocal of its level frequency $f_{d,j}$, in a way that the more the feature $j$ is used at level $d$ in the $TE$, the less is the penalization for using it in the optimal tree. We restrict this penalization only to the most used features per level, namely for $j\in\J_\gamma(d)$.
The penalization term is

\begin{align}
V(T)= \sum_{d\in\D}  \sum_{j\in \J_\gamma(d)}\frac{1}{f_{d,j}} \sum_{t\in\Tb(d)}s_{t,j}. \nonumber
\end{align}

Thus, the overall objective function is:

$$\min \ds \frac 1 2 \sum_{i\in\I} p^{i}\hat{y}^{i}\Big(\hat{y}^{i} - \sum_{\ell\in\Tl} c_{\ell}z_{i,\ell}\Big) + \alpha \sum_{d\in\D}\sum_{j\in \J_\gamma(d)} \frac{1}{f_{d,j}} \sum_{t\in{\Tb}(d)} s_{t,j} ,$$

where $\alpha$ is a penalty hyperparameter to control the trade-off between the two objectives.

\subsection{The basic MILP model of \texttt{MIRET}}

The MILP formulation of the basic-\texttt{MIRET} is the following:
{\linespread{1.5}\selectfont
\footnotesize
\begin{alignat*}{4}
    (\text{\texttt{b-MIRET}}) \quad
      \min \limits_{{a,b,z,s} } \quad & \ds \frac 1 2 \sum_{i\in\I} p^{i}\hat{y}^{i}\Big(\hat{y}^{i} - \sum_{\ell\in\Tl} c_{\ell}z_{i,\ell}\Big) + \alpha\sum_{d\in\D}\sum_{j\in \J_\gamma(d)} \frac{1}{f_{d,j}} \sum_{t\in{\Tb}(d)} s_{t,j} \nonumber\\ 
     \text{s.t.} \quad
& a_t^T x^i + b_t - \varepsilon \geq -(|J|+1+\varepsilon) \Big(1 - \sum_{\ell\in\S_R(t)} z_{i,\ell}\Big) && \forall t \in \Tb, &&\quad\forall i\in\I,\\
     & a_t^T x^i + b_t \leq (|J|+1)\Big(1 - \sum_{\ell\in\S_L(t)} z_{i,\ell}\Big) &&\forall t \in \Tb ,&&\quad \forall i\in\I, \\
    & \ds\sum_{\ell\in\Tl} z_{i,\ell} = 1  && \forall i\in\I, \\
   & z_{i,\ell} = z_{k,\ell} && \forall \ell\in\Tl, &&\quad  \forall i,k \in \U_{TE}{(\overline{m}_{TE})}, \\                                                                                                                
    & a_{t,j} = 0,  &&\forall d\in\mathcal{D},   &&\quad \forall t\in {\Tb}(d),\quad \forall j\in \J\setminus \J_\gamma(d),\\
    & s_{t,j} = 0  &&\forall d\in\mathcal{D},&&\quad \forall t\in {\Tb}(d),\quad \forall j\in \J\setminus \J_\gamma(d), \\
    & -s_{t,j} \leq a_{t,j} \leq s_{t,j}  && \forall t\in\Tb, &&\quad \forall j\in\J, \\
      &-1\le a_{tj}\le 1 &&   \forall t \in \Tb\ &&\quad \forall j\in\J,\\
    &-1\le b_t\le 1 &&   \forall t \in \Tb,\\
     & z_{i,\ell} \in\{0,1\} &&\forall \ell \in \Tl,  &&\quad \forall i\in\I,\\
    & s_{t,j} \in\{0,1\} && \forall t\in\Tb, &&\quad  \forall  j\in\J.
\end{alignat*}
}
Solving this model gives an optimal tree $T^*$ which for all nodes $t\in \B$ at a node $t$ has the splits
$$
\text{if } {a^*_t}^T x+ b^*_t  \begin{cases}
    \leq 0 \cr
 >0            
\end{cases},$$
Following the path from the root to leaves gives the decision rules that lead to a classification of a sample.
In the next section, we slightly modify the formulation in order to gain on the performance of out-of-shell MILP solvers.

\subsection{\lnote{Comparison between other state-of-the-art formulations for multivariate OCTs}}

\lnote{In this section, we highlight some technical differences between our \texttt{MIRET} forumlation and other recent MIP formulations for multivariate trees in the literature. In particular, we address the comparison with the OCT-H model presented in \cite{Bertsimas2017OptimalClassification} and the S-OCT formulation proposed in \cite{boutilier2022shattering}. 
The goal of our model is to create a single multivariate surrogate tree that serves as an interpretable representer of a more complex tree ensemble classifier. 
To this aim, we embed specific TE metrics into our model, that are specially designed with the aim of creating 'born-again' trees.
In contrast,   OCT-H and S-OCT  aim to provide a standalone predictor. 
In principle, TE metrics that have been used in the \texttt{MIRET} model  (maximization of the fidelity, proximity and frequency constraints, features frequency and class probability driven penalization) could be integrated into any mixed-integer model for multivariate optimal trees. 
For this reason, the comparison below focus solely on technical aspects related to generic tree construction (objective function and routing constraints) excluding 
 the TE-based metrics that are a novelty specific to our formulation, but can be applied beyond our specific model.}

\lnote{Concerning routing constraints, they involve binary variables, often leading to large and complex mixed-integer models. 
Aiming to reduce the model complexity to obtain a more tractable problem,
we mainly follow the modeling approaches for routing and assignment constraints used in \cite{MarginOptimalClassificationTrees},  
adapting them to suit our specific tree model.
{We opted for this approach because it allowed us to reduce the number of binary variables and associated constraints.}
 \lnote{
Indeed, as in OCT-H \cite{Bertsimas2017OptimalClassification}, we defined $z$ binary variables solely over the leaf nodes, resulting in $|\I|\cdot|\Tl|=|\I|\cdot2^D$ variables. 
On the other hand, in S-OCT, such binary variables are defined \lnote{over branch and leaf nodes} totaling \lnote{$|\I|\cdot(|\Tb|+|\Tl|)=|\I|\cdot(2^{D+1}-1)$.}}
{Second, as in \cite{MarginOptimalClassificationTrees}, we model routing constraints by using the \lnote{sets of left and right sub-leaves} $\S_R(t), \S_L(t), t\in\Tb$, generating a total of  {$2\cdot|\I|\cdot|\B|=|\I|\cdot(2^{D+1}-2)$} constraints.
In OCT-H instead, routing constraints are defined for each leaf considering its set of right and left ancestors leading to a higher number of constraints ($|\I|\cdot|\Tl|\cdot D=|\I|\cdot2^D\cdot D$) compared to \texttt{MIRET} formulation.
In constrast, in S-OCT, routing constraints are modeled for each branch node amounting to $2\cdot |\I|\cdot|\Tb|=|\I|\cdot(2^{D+1}-2)$. Additionally, other $|\I|\cdot|\Tb|=|\I|\cdot(2^D-1)$ constraints are needed to ensure that each sample assigned to a branch node $t$ is routed to either the left or right child $t$. In Table \ref{tab: dimension}, a summary of the number of integer variables and routing constraints is presented.}\\
\gnote{Finally, the objective function in optimal tree model balances a trade-off between the empirical misclassification loss and a regularization term usually expressed as a norm of the weights of the hyperplane splits in the tree. In \texttt{MIRET}, the second term addresses the local sparsity of the
tree model by penalizing the $\ell_1$-norm of the weights and uses a linearization using binary variables $s$, similar to OCT-H.
However, in  \texttt{MIRET}, each term of the linearization is weighted in the objective function by the TE-based frequency \eqref{eq: miret_frequency} of the corresponding feature.
In OCT-H, the regularization term does not consider the TE tailored weights. 
On the other hand, in S-OCT, the second term addresses the complexity of the tree, i.e. the number of effective splits in the tree, with the use of binary variables that model whether a node applies a split or not.}}

\lnote
{As a further check to understand the role played by the TE metrics in the construction of the represented tree, we present in \ref{app: appendix1} a comparison among the \texttt{MIRET} tree and a simple MIP model excluding TE metrics. The results show that TE metrics help in defining a more suitable representer tree and in facilitating the solution.}

\begin{center}
\begin{table}[ht]
\renewcommand\arraystretch{1.2}
\centering
\begin{tabular}{lllllll}
\multicolumn{1}{c}{Class}      & \multicolumn{1}{c}{} & \multicolumn{5}{c}{Size}                                                                                          \\ \hline
\multicolumn{1}{c}{}           & \multicolumn{1}{c}{} & \multicolumn{1}{c}{\texttt{MIRET}} &  & \multicolumn{1}{c}{OCT-H} &  & \multicolumn{1}{c}{S-OCT} \\ \hline\hline
Integer variables (assignment) &                      & $|\I|2^D$                            &  & $|\I|2^D$                            &  & $|\I|(2^{D+1}-1)$                        \\
Routing constraints            &                      & $|\I|(2^{D+1}-2)$                    &  & $|\I|2^DD$                           &  & $|\I|(2^{D+1}-2)+|\I|(2^D-1)$        \\ \hline
\end{tabular}\caption{\lnote{Number of integer variables (modeling assignment of samples) and routing constraints in  \texttt{MIRET}, OCT-H and S-OCT respectively.}}
    \label{tab: dimension}
\end{table}
\end{center}


\section{Improving the basic \texttt{MIRET} formulation} \label{sec:strenghen}

\subsection{Existence of at least one split}
As a first step, we aim to avoid the existence of a feasible solution that does not provide any split, namely a dummy tree. Indeed, we observe that there exists a feasible (obviously not optimal) solution with  $s_{j,t}=0$ for all $j\in\J$ and $t\in\B$ (so that the penalty term $V(T)=0$) and $a_t=0$ for all $t\in \B$ where $b_t, z_{i,\ell}$ can be fixed to any feasible value. All these feasible solutions do not apply a partition of the samples, thus resulting in a dummy tree. 
In particular, this also implies that when using a continuous linear relaxation to compute a lower bound, we can obtain a fractional feasible solution 
$z_{i,\ell}\in [0,1]$  such that the assignment constraints are satisfied and 
$\hat{y}^{i} = \sum_{\ell\in\Tl} c_{\ell}z_{i,\ell}$ which is of course optimal for the linear relaxation since the value of the objective function is zero.
This gives a trivial lower bound.

In order to avoid such a trivial solution, we add the constraint
\begin{align}\label{eq: cons trivial}
\sum_{t\in\Tb}\sum_{j\in \J}  s_{t,j}\ge 1 \nonumber
\end{align}
which forces at least one branch node in the tree to use at least one feature (the tree has at least one univariate split).
Of course, we could enforce a more stringent constraint by increasing the right-hand side to a value greater than $1$ which is usually the case in optimal tree solutions. 

\subsection{Breaking symmetries}



In a MILP model, symmetries refer to the existence of multiple solutions that have the same objective value but differ only in the values assigned to a set of integer variables \cite{Margot2010}.
The presence of symmetries can make it difficult to solve the problem, as they can lead to a large number of equivalent solutions, and to a huge branching tree to be explored by the optimization algorithm.

Our basic formulation is affected by this issue in a specific way.
Indeed, our classification tree is symmetric in the sense that the same value of the objective function is obtained for every solution that permutes components of $z$ variables in such a way that
\begin{enumerate}[(i)]
\item the final partition of samples, and so the final predictions, remains unchanged and,
\item  the same $s$ variables activate, i.e. the same features are used at each level $d$ of the tree.
\end{enumerate}


Indeed, in any branch node $t\in\Tbf$, samples assigned to its right and left sub-trees can be swapped using the same setting of the variables $s$ obtaining equivalent solutions.


However, this is no longer true for nodes at the branching level adjacent to the leaves $\Tbl$ due to the pre-assigned labels $c$ to leaves. Hence, for each 
 node $t\in\Tbl$, the subset of samples $I_t$ is forced to follow either the left or the right branch of the node  by the leaves labels. 




In order to avoid the presence of multiple equivalent solutions obtained by swapping branches,  we require 
the following symmetry-breaking constraints:
\begin{align}
& b_t\ge 0, & \forall t \in \Tbf \nonumber
\end{align}
This set of constraints enforces a non-negative intercept for the splitting hyperplanes at node $t\in \Tbf$ thus avoiding the presence of symmetric solutions obtainable by reversing the splitting rules.
In this way, we reduce the set of feasible solutions by removing redundant solutions which produce the same final partitions of samples using the same features in the tree.

\subsection{Adding branching rules}


Exploiting the nested tree structure and the routing constraints \eqref{cons: routing left} \eqref{cons: routing right}, we introduce for each $t\in\Tbf$ and each sample $i\in\I$ two additional boolean variables $q^i_L(t), q^i_R(t) \in \{0,1\}$ with the constraints
$$
 q^i_L(t) - \sum_{\ell\in\S_L(t)} z_{i,\ell}=0 ;\qquad\qquad
 q^i_R(t) - \sum_{\ell\in\S_R(t)} z_{i,\ell}=0.
$$

Hence, we have
$$
    q^i_{L/R}(t) =
    \begin{cases} 
      1 & \text{if sample $i$ is assigned to a leaf $\ell\in \S_{L/R}(t)$} \\
      0 & \text{otherwise} 
   \end{cases},
$$
and routing constraints \eqref{cons: routing right}-\eqref{cons: routing left}  become
\begin{align*}
     & a_t^T x^i + b_t \leq M_{\H_L}\Big(1 -  q^i_L(t)\Big), && \forall t \in \Tb , &&\quad \forall i\in\I\\
 & a_t^T x^i + b_t - \varepsilon \geq -M_{\H_R} \Big(1 -  q^i_R(t)\Big), && \forall t \in \Tb, &&\quad \forall i\in\I
 \end{align*}
and 
the assignment constraints \eqref{cons: fin assignment} are written as
\begin{align*}
&  q^i_L(0)+q^i_R(0) = 1,  &&  \forall i\in\I.
\end{align*}




In the branching procedure,  fixing a variable $q_L^i(t)=0$ (or  $q_R^i(t)=0$)  forces all the  $z_{i,\ell}=0$ with $\ell\in \S_L(t)\ (\text{or }S_R(t))$, thus avoiding the need to explore branching on the single variables $z_{i,\ell}$.

Although the large increase in variables and  constraints, the numerical experiments that we performed seem to confirm the effectiveness of having such variables.

We also add parenting constraints from node $t\in\B(d)$ to its children $2t+1, 2t+2 \in \B(d+1)$ that connects the $q_L$ and $q_R$ variables along the subtree rooted at $t$. These are $|\I|\cdot|\Tbf| $ conditions expressed as follow

\begin{align*}
& q^i_L(t) = q^i_L(2t+1) +  q^i_R(2t+1), &&\forall t\in\Tbf, &&\forall i\in\I,\\
 & q^i_R(t) = q^i_L(2t+2) +  q^i_R(2t+2), &&\forall t\in\Tbf, &&\forall i\in\I.
\end{align*}

in Table \ref{tab:variables_inmproved} we report the additional variables.

\begin{table}[ht!]
\renewcommand\arraystretch{1.4}
\footnotesize
\centering
\begin{tabular}{lll}
\textbf{Variable}                & & \textbf{Description}                                             \\ \hline\hline
$q^i_L(t) \in \{0,1\}$    & & if sample $i$ is assigned to a leaf $\ell\in \S_{L}(t)$  
\\
$q^i_R(t) \in \{0,1\}$    & & if sample $i$ is assigned to a leaf $\ell\in \S_{R}(t)$    
\\
\end{tabular}\caption{Additional decision variables in improved {\texttt{MIRET}}}
    \label{tab:variables_inmproved}
\end{table}

The improved MILP formulation \texttt{MIRET} is the following:
{\linespread{1.5}\selectfont
\footnotesize
\begin{alignat*}{4}
    (\text{\texttt{MIRET}}) \quad
      \min \limits_{{a,b,z,s,q_L,q_R} } \quad & \ds \frac 1 2 \sum_{i\in\I} p^{i}\hat{y}^{i}\Big(\hat{y}^{i} - \sum_{\ell\in\Tl} c_{\ell}z_{i,\ell}\Big) + \alpha\sum_{d\in\D}\sum_{j\in \J_{\gamma}(d)} \frac{1}{f_{d,j}} \sum_{t\in{\Tb}(d)} s_{t,j} \nonumber\\ 
     \text{s.t.} \quad
 & a_t^T x^i + b_t - \varepsilon \geq -(|J|+1+\varepsilon) \Big(1 -  q^i_R(t)\Big) && \forall t \in \Tb, && \quad \forall i\in\I,\\
& a_t^T x^i + b_t \leq (|J|+1)\Big(1 -  q^i_L(t)\Big) && \forall t \in \Tb , &&\quad \forall i\in\I,\\
&  q^i_L(0)+q^i_R(0) = 1  &&  \forall i\in\I, \\
  & q^i_L(t) = q^i_L(2t+1) +  q^i_R(2t+1), &&\forall t\in\Tbf, &&\quad\forall i\in\I,\\
 & q^i_R(t) = q^i_L(2t+2) +  q^i_R(2t+2), &&\forall t\in\Tbf, &&\quad\forall i\in\I, \\
    & q^i_L(t) = \sum_{\ell\in\S_L(t)} z_{i,\ell} &&\forall t \in \Tb ,&&\quad \forall i\in\I, \\
    & q^i_R(t) = \sum_{\ell\in\S_R(t)} z_{i,\ell} &&\forall t \in \Tb ,&&\quad \forall i\in\I, \\
       &\sum_{t\in\Tb}\sum_{j\in \J}  s_{t,j}\ge 1,\\
       & b_t\ge 0, && \forall t \in \Tbf,\\
   & z_{i,\ell} = z_{k,\ell} && \forall \ell\in\Tl, &&\quad \forall i,k \in \U_{TE}(\overline{m}_{TE}), \\
    & a_{t,j} = 0,  &&\forall d\in\mathcal{D}, &&\quad\forall t\in {\Tb}(d), \quad\forall j\in \J\setminus \J_\gamma(d),\\
    & s_{t,j} = 0  &&\forall d\in\mathcal{D}, &&\quad\forall t\in {\Tb}(d), \quad \forall j\in \J\setminus \J_\gamma(d),\\
        & -s_{t,j} \leq a_{t,j} \leq s_{t,j},  && \forall t\in\Tb, &&\quad \forall j\in\J, \\
      &-1\le a_{t,j}\le 1 &&   \forall t \in \Tb\ &&\quad \forall j\in\J,\\
    &-1\le b_t\le 1 &&   \forall t \in \Tb,\\
     & z_{i,\ell} \in\{0,1\} &&\forall \ell \in \Tl,  &&\quad \forall i\in\I,\\
    & s_{t,j} \in\{0,1\} && \forall j\in\J, &&\quad  \forall t\in\Tb.
\end{alignat*}
}

\section{Computational experience}\label{sec:numerical}

In this section, we present different computational results in order to evaluate the performances of the approach proposed.
We use a Random Forest as $TE$, and in particular, we use the RF method as implemented in \texttt{sklearn.ensemble.RandomForestClassifier} \cite{scikit-learn} with the following setting:
\begin{itemize}
    \item maximum depth $D\in\{2,3,4\};$
\item $|\E|=100;$
\item random sampling {of the features} deactivated so that weights $w^e=1$, $e\in\E$ in  \eqref{level_frequency} and \eqref{eq:proximity}.
\end{itemize} 

We denote the target forest as $\widehat{RF}$. 
The mixed-integer programming model \texttt{MIRET} is coded in Python and ran 
on a server Intel(R) Xeon(R) Gold 6252N CPU processor at 2.30 GHz and 96 GB of RAM.
The MILP is solved using Gurobi 10.0.0 with standard settings. We set a time limit {of 1 hour} for each model optimization.

\subsection{Datasets}

We selected 10 datasets from UCI Machine Learning repository \cite{UCI2019} related to binary classification tasks, and we normalized the feature values of each dataset in the interval 0-1. Information about the datasets considered is reported in Table \ref{table: datasets}.
Each dataset was split into training ($80\%$) and test ($20\%$) sets.


\begin{table}[ht!]
\centering
\scalebox{1}{
\renewcommand\arraystretch{1.3}
\begin{tabular}{lcccc}
Dataset     &  & $|\I|$ & $|\J|$ & Class (\%) \\ \hline
Cleveland   &  & 297    & 13  & 53.9/46.1  \\
Diabetes    &  & 768    & 8   & 65.1/34.9  \\
German      &  & 1000   & 20  & 30/70      \\
Heart       &  & 270    & 13  & 55.6/44.4  \\
IndianLiver &  & 579    & 10  & 71.5/28.5  \\
Ionosphere  &  & 351    & 34  & 35.9/64.1  \\
Parkinson   &  & 195    & 22  & 24.6/75.4  \\
Sonar       &  & 208    & 60  & 53.4/46.6  \\
Wholesale   &  & 440    & 7   & 32.3/67.7  \\
Wisconsin   &  & 569    & 30  & 37.3/62.7 
\end{tabular}}
\caption{Characteristics of the datasets.}
\label{table: datasets}
\end{table}

In order to have a glimpse view of the role of the $\widehat{RF}$ information that we used in the model, we calculate  the distribution of the proximity measures and of the probability class on the training set, and we plot them using violin plots \cite{hintze1998violin} in  Figures \ref{fig:Prox_meas_dist} and \ref{fig:Prob_meas_dist} respectively.
As regard the level frequencies, we also calculate them for the \gnote{$\widehat{RF}$}, and we report them together with those used in the optimal tree obtained with \texttt{MIRET} in Figure \ref{fig:Feat_selected}. 

\begin{figure}
\centering
\includegraphics[width=1.\textwidth]{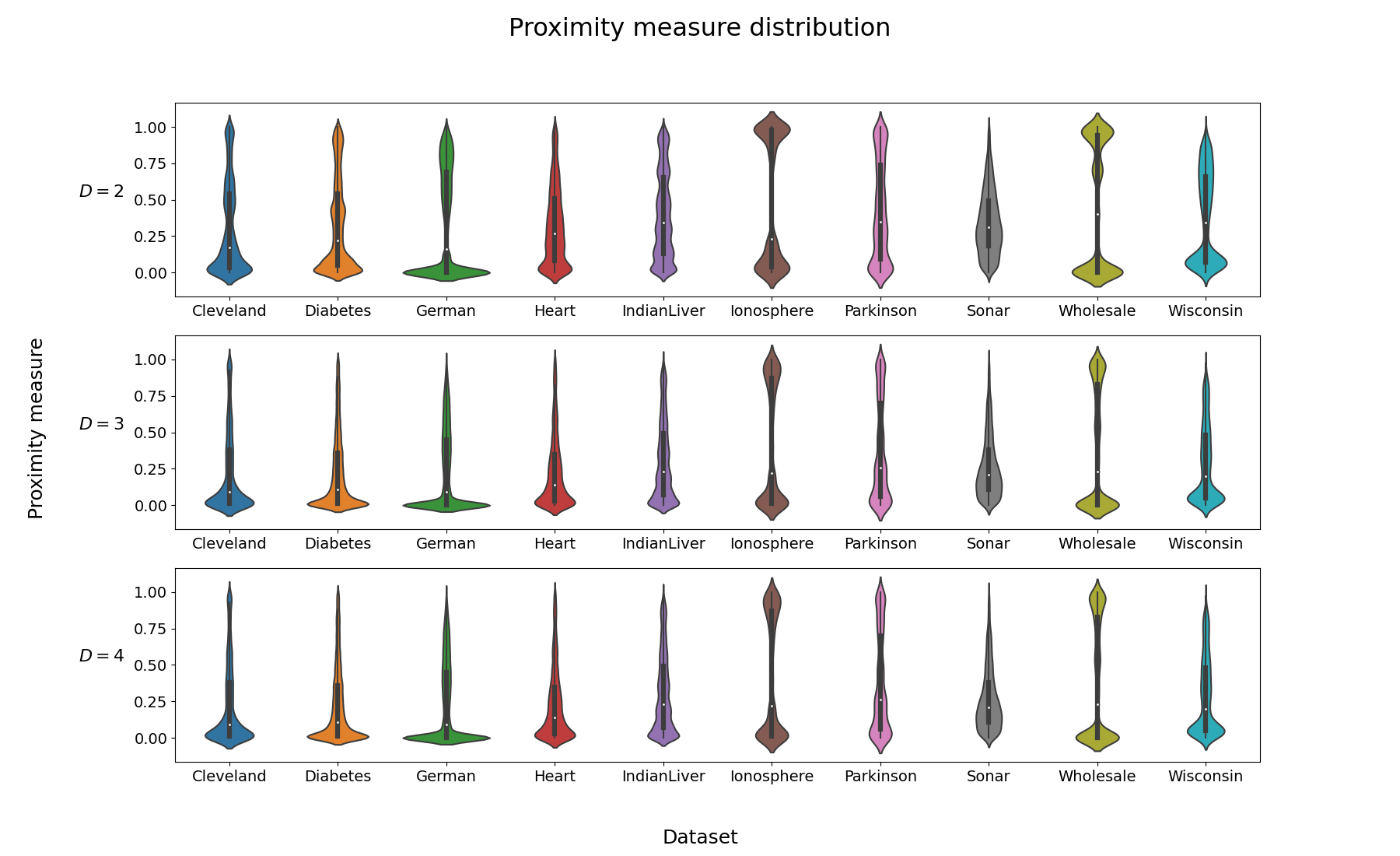}
\caption{Distribution of proximity measures of pair of training set samples}
\label{fig:Prox_meas_dist}
\end{figure}

\begin{figure}
\centering
\includegraphics[width=1.\textwidth]{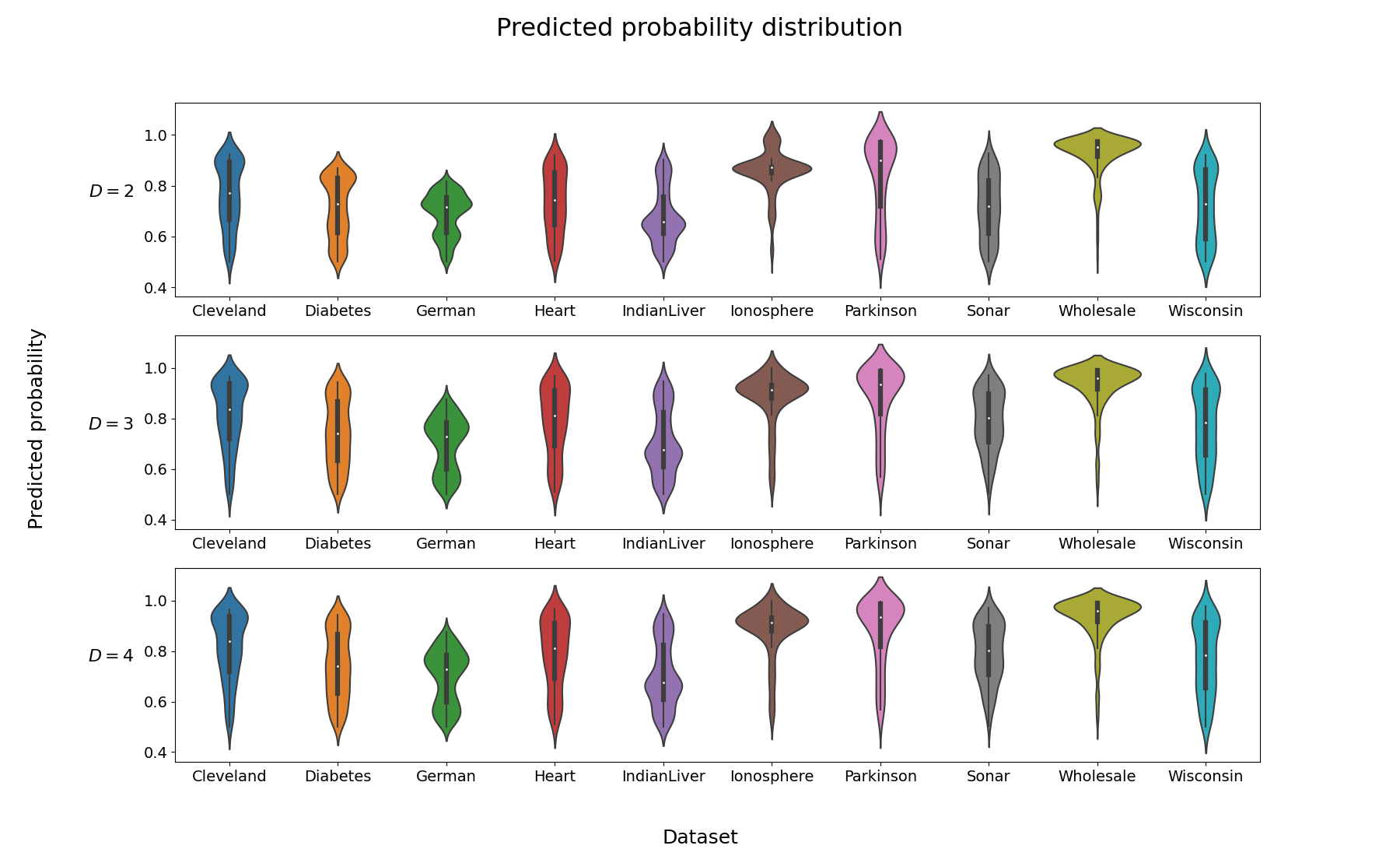}
\caption{Distribution of predicted probabilities of training set samples}
\label{fig:Prob_meas_dist}
\end{figure}

\subsection{Hyperparameters setting}\label{sub: hyperparam_set}

For hyperparameters $\gamma_d$, $d\in\D$, \mnote{$\overline{m}_{TE}$} and $\alpha$, we performed a tailored  tuning of the model using a grid search
within a $k$-fold cross-validation \gnote{with $k=4$}. 

In particular, the grid values on the thresholds $\gamma_d$ in \eqref{def:gamma_d} are defined considering a percentage of the frequencies values of features used in the level $d$ of the $\widehat{RF}$. Let  $\F(d)^+=\{f_{d,j}: f_{d,j}>0, \ j\in\J\}$  be the set of positive level frequencies, then $\gamma_d$, $d\in\D$ are computed as the $h$-th percentile of $\F(d)^+$. Further, we also include the value $\gamma_d = 0$, to consider the case when \texttt{MIRET} can use all the features used in $\widehat{RF}$, namely any feature in $\F(d)^+$.

We perform a grid search using these values:

\begin{itemize}
\item penalty parameter $\alpha \in \{0.2,0.4,0.5,0.6,0.8\}$;
\item \lnote{$\gamma_d=0$ and $\gamma_d$  as the $h$-th percentile of $\F(d)^+$ with $h\in\{100/2,100/3,100/4\}$;
\item proximity threshold $\overline{m}_{\widehat{RF}} \in \{0.85,0.90,1.00\}$.}

\end{itemize}
We evaluated the average validation fidelity and the average sparsity of the trained model.
In the end, we selected hyperparameters that provide the best balance between the two objectives.
The final hyperparameters setting used in the computational results are reported in the Appendix \ref{app: appendix2} in Table \ref{tab: hyperparameters}.


\subsection{Results} \label{sec: results}

We analyze both the optimization performances, namely the viability and efficiency of using a MILP formulation for finding an optimal representative tree $T^*$ of the random forest $\widehat{RF}$, and the predictive performances. 
We select depths up to $D=4$. Indeed, aiming to obtain an interpretable tree model and accounting for the fact that we can have multivariate splittings, depth $D\in\{2,3\}$ are the most significant to get real interpretability. However, to check scalability in the solution of the model, we also use $D=4$.

 Once an optimal tree $T^*$  is obtained for $D\in\{2,3,4\}$ solving the MILP problem \texttt{MIRET}, following \cite{guidotti2018survey}, we evaluate the performances 
 both with respect to the $\widehat{RF}$ predictions (fidelity), and with respect to the ground truth labels (accuracy) both on the training and on the test sets. 
\begin{itemize}
    \item  Fidelity is measured as 
    $$\text{FID}_{\widehat{RF}}=1 -  \frac{1}{2 |\I|} \displaystyle\sum_{i\in\I}  F_{\widehat{RF}}(x^i)\Big(F_{\widehat{RF}}(x^i) - {F_{T^*}(x^i)}\Big);$$
    \item \texttt{MIRET} accuracy with respect to the ground truth $y^i$, calculated  as 
    $$\text{ACC}_{\texttt{MIRET}}=1- \frac{1}{2|\I|}\displaystyle\sum_{i\in\I}  y^i\Big(y^i - {F_{T^*}(x^i)}\Big).$$
 \end{itemize}
We also evaluated as a term of comparison the $\widehat{RF}$ accuracy with respect to the ground truth $y^i$, calculated as 
$$\text{ACC}_{\widehat{RF}} =1 - \frac{1}{2|\I|}\displaystyle\sum_{i\in\I}  y^i\Big(y^i - F_{\widehat{RF}}(x^i)\Big) .$$

\begin{table}[ht!]
\centering
\scalebox{0.65}{
\renewcommand\arraystretch{1.5}
\begin{tabular}{lcccccccccccccc}
              & \multicolumn{4}{c}{$D=2$}                                                     &  & \multicolumn{4}{c}{$D=3$}                                                     &  & \multicolumn{4}{c}{$D=4$}                                                     \\ \hline
              & \multicolumn{2}{c}{Time}              & \multicolumn{2}{c}{Gap}               &  & \multicolumn{2}{c}{Time}              & \multicolumn{2}{c}{Gap}               &  & \multicolumn{2}{c}{Time}              & \multicolumn{2}{c}{Gap}               \\
Dataset       & $\texttt{b-MIRET}$ & $\texttt{MIRET}$ & $\texttt{b-MIRET}$ & $\texttt{MIRET}$ &  & $\texttt{b-MIRET}$ & $\texttt{MIRET}$ & $\texttt{b-MIRET}$ & $\texttt{MIRET}$ &  & $\texttt{b-MIRET}$ & $\texttt{MIRET}$ & $\texttt{b-MIRET}$ & $\texttt{MIRET}$ \\ \hline
Cleveland     & \textbf{28.4}      & 44.6             & \textbf{0.0 }               & \textbf{0.0}              &  & \underline{3600}   & \textbf{1942.2}  & 43.6               & \textbf{0.0}     &  & \underline{3600}   & \underline{3600} & {89.8}             & {\textbf{71.5}}  \\
Diabetes      & \underline{3600}   & \textbf{1373.3}  & 39.5               & \textbf{0.0}     &  & \underline{3600}   & \underline{3600} & 100.0              & \textbf{96.1}    &  & \underline{3600}   & \underline{3600} & 100.0              & {\textbf{98.6}}  \\
German        & \textbf{11.8}      & 12.3             & \textbf{0.0 }               & \textbf{0.0}              &  & \underline{3600}   & \underline{3600} & 99.6               & \textbf{73.6}    &  & \underline{3600}   & \underline{3600} & {94.6}             & {\textbf{93.4}}  \\
Heart         & 160.0              & \textbf{0.7}    & \textbf{{0.0}}              & \textbf{0.0}              &  & \underline{3600}   & \underline{3600} & 87.6               & \textbf{73.2}    &  & \underline{3600}   & \underline{3600} & {\textbf{87.4}}    & {95.0}           \\
IndianLiver   & \underline{3600}   & \underline{3600} & 80.5               & \textbf{75.6}    &  & \underline{3600}   & \underline{3600} & 83.2               & \textbf{77.1}    &  & \underline{3600}   & \underline{3600} & 100.0              & {\textbf{98.3}}  \\
Ionosphere    & 4.4                & \textbf{3.4}     & \textbf{0.0}                & \textbf{0.0}              &  & \textbf{107.4}     & 204.4            & \textbf{0.0}                & \textbf{0.0}              &  & \underline{3600}   & \textbf{947.6}   & {56.6}             & {\textbf{0.0}}   \\
Parkinson     & 5.6                & \textbf{4.6}     & \textbf{0.0}                & \textbf{0.0}              &  & \underline{3600}   & \underline{3600} & 58.2               & \textbf{33.8}    &  & \underline{3600}   & \underline{3600} & {\textbf{92.6}}    & {96.3}           \\
Sonar         & \underline{3600}   & \underline{3600}  & 38.6               & \textbf{12.0}     &  & \underline{3600}   & \underline{3600} & \textbf{89.3}      & 93.2             &  & \underline{3600}   & \underline{3600} & {97.9}             & {\textbf{97.1}}  \\
Wholesale     & 1.4                & \textbf{0.7}     &\textbf{ 0.0  }              & \textbf{0.0}             &  & 70.8               & \textbf{18.0}    & \textbf{0.0}                &\textbf{ 0.0  }            &  & \underline{3600}   & \textbf{554.2}   & {59.5}             & {\textbf{0.0}}   \\
Wisconsin     & 61.3               & \textbf{8.2}     & \textbf{0.0}                &\textbf{0.0}              &  & 814.2              & \textbf{385.5}   & \textbf{0.0}                & \textbf{0.0}     &  & \underline{3600}   & \underline{3600} & 99.5               & {\textbf{98.1}}  \\ \hline
\textbf{Mean} & 1107.3             & \textbf{864.8}   & 22.7               & \textbf{8.8}     &  & 3163.4             & \textbf{2415.0}  & 56.2               & \textbf{44.7}    &  & 3600.0             & \textbf{3030.2}  & 87.8               & \textbf{74.8}    \\ \hline
\end{tabular}}
\caption{\lnote{\texttt{b-MIRET} and \texttt{MIRET} comparison on
computational times (s) and MIP Gap values (\%); in boldface the winning values.}}
\label{table: comparison_form}
\end{table}

In Table \ref{table: comparison_form}, we compare the optimization performances of the basic \texttt{MIRET} model and the improved \texttt{MIRET} version.
 For each dataset and each $D$ (a total of 30 MILP problems), we report the computational time (s) and the optimality gap for solving the corresponding MILP.
 The returned solution $T^*$  is not always certified as the global optimum of the problem. 
As expected, problem hardness increases with the size of the problem. 
However, the improved version of \texttt{MIRET} has better performance because it closes the gap on \gnote{four} additional problems, improves the gap on \lnote{85.0\% of  the problems not closed by \texttt{b-MIRET}}, and improves the time on about \lnote{70.0\%} of the closed ones. 
In particular, for $D=2$, we obtain gap zero on \lnote{eight} out of the ten datasets; for $D=3$ optimality is certified  on {four} of the ten datasets and finally only on \lnote{two} problems for $D=4$. 
Despite the gap is not zero, the quality of the solutions in terms of predictive performance on the training set is outstanding, as reported in Table \ref{table: train_acc}.

\begin{table}[ht!]
\centering
\scalebox{0.85}{
\renewcommand\arraystretch{1.3}
\begin{tabular}{lcccccccccccc}
            &  & \multicolumn{3}{c}{$D = 2$}                            &  & \multicolumn{3}{c}{$D = 3$}                            &  & \multicolumn{3}{c}{$D = 4$}                            \\ \hline
            &  &                    & \multicolumn{2}{c}{ACCURACY}      &  &                    & \multicolumn{2}{c}{ACCURACY}      &  &                    & \multicolumn{2}{c}{ACCURACY}      \\
Dataset     &  & FID-$\widehat{RF}$ & $\texttt{MIRET}$ & $\widehat{RF}$ &  & FID-$\widehat{RF}$ & $\texttt{MIRET}$ & $\widehat{RF}$ &  & FID-$\widehat{RF}$ & $\texttt{MIRET}$ & $\widehat{RF}$ \\ \hline
Cleveland   &  & 96.6               & 81.9             & 81.9           &  & 97.9               & 85.7             & 87.8           &  & 93.7               & 85.7             & 92.0           \\
Diabetes    &  & 91.9               & 73.5             & 75.4           &  & 95.4               & 75.9             & 79.2           &  & 92.3               & 76.9             & 83.2           \\
German      &  & 98.0               & 66.8             & 66.0           &  & 95.4               & 68.5             & 71.6           &  & 87.8               & 68.8             & 78.0           \\
Heart       &  & 77.3       & 74.5     & 84.3           &  & 95.4               & 86.6             & 89.4           &  & 93.1               & 87.0             & 94.0           \\
IndianLiver &  & 97.2               & 67.6             & 67.8           &  & 95.7       & 67.2     & 70.6           &  & 93.5               & 72.1             & 77.3           \\
Ionosphere  &  & 98.6               & 89.6             & 91.1           &  & 95.0               & 89.6             & 94.6           &  & 95.0               & 92.5             & 96.8           \\
Parkinson   &  & 93.6               & 84.6             & 88.5           &  & 89.1               & 85.9             & 95.5           &  & 96.2               & 94.9             & 98.7           \\
Sonar       &  & 86.7       & 79.5     & 90.4           &  & {83.3}             & {83.3}           & 100.0          &  & 81.9               & 81.9             & 100.0          \\
Wholesale   &  & 99.7               & 92.6             & 92.3           &  & 99.7               & 92.6             & 92.9           &  & 97.2               & 92.3             & 95.2           \\
Wisconsin   &  & 99.1               & 95.8             & 96.7           &  & 96.9               & 95.6             & 98.7           &  & 98.5               & 97.8             & 99.3           \\ \hline
\end{tabular}}
\caption{\mnote{Predictive performances on the training set. FID-$\widehat{RF}$: Fidelity with respect to $\widehat{RF}$ (\%); \texttt{MIRET} ACCURACY  and  $\widehat{RF}$ ACCURACY with respect to the ground truth (\%).}}
\label{table: train_acc}
\end{table}

In Table \ref{table: test_acc} we report the out-of-sample predictive performances. The average fidelity of \texttt{MIRET} (FID-$\widehat{RF}$) is quite high, meaning that it fairly mimics the $\widehat{RF}$ classification, with an average of \gnote{94.1\%} at depth 2, \gnote{94.1\%}  at depth 3, and \gnote{92.4\%} at depth 4.

Additionally, our model has retained on average the generalization capabilities of $\widehat{RF}$.
As a matter of fact, it can be easily seen that the accuracy of \texttt{MIRET} is similar to the accuracy of $\widehat{RF}$ along all the datasets. 
Thus, our model gains in interpretability as it provides a single optimal tree with more straightforward decision paths while being capable of fairly closely reproducing the predictive performances of the Random Forest.

\begin{table}[ht!]
\centering
\scalebox{0.85}{
\renewcommand\arraystretch{1.3}
\begin{tabular}{lcccccccccccc}
            &  & \multicolumn{3}{c}{$D = 2$}                            &  & \multicolumn{3}{c}{$D = 3$}                            &  & \multicolumn{3}{c}{$D = 4$}                            \\ \hline
            &  &                    & \multicolumn{2}{c}{ACCURACY}      &  &                    & \multicolumn{2}{c}{ACCURACY}      &  &                    & \multicolumn{2}{c}{ACCURACY}      \\
Dataset     &  & FID-$\widehat{RF}$ & $\texttt{MIRET}$ & $\widehat{RF}$ &  & FID-$\widehat{RF}$ & $\texttt{MIRET}$ & $\widehat{RF}$ &  & FID-$\widehat{RF}$ & $\texttt{MIRET}$ & $\widehat{RF}$ \\ \hline
Cleveland   &  & 98.3               & 78.3             & 80.0           &  & 95.0               & 76.7             & 81.7           &  & 93.3               & 76.7             & 83.3           \\
Diabetes    &  & 89.6               & 75.3             & 68.8           &  & 95.5               & 77.3             & 76.6           &  & 94.2               & 78.6             & 76.6           \\
German      &  & 99.0               & 65.5             & 65.5           &  & 97.0               & 68.5             & 67.5           &  & 87.0               & 68.0             & 72.0           \\
Heart       &  & 77.8       & 77.8     & 85.2           &  & 90.7               & 79.6             & 85.2           &  & 90.7               & 79.6             & 85.2           \\
IndianLiver &  & 94.8               & 65.5             & 65.5           &  & 94.0       & 64.7     & 65.5           &  & 96.6               & 63.8             & 63.8           \\
Ionosphere  &  & 100.0              & 94.4             & 94.4           &  & 95.8               & 94.4             & 93.0           &  & 87.3               & 88.7             & 93.0           \\
Parkinson   &  & 97.4               & 76.9             & 79.5           &  & 100.0              & 79.5             & 79.5           &  & 97.4               & 84.6             & 82.1           \\
Sonar       &  & 90.5       & 71.4     & 81.0           &  & {83.3}             & {83.3}           & 81.0           &  & 88.1               & 71.4             & 78.6           \\
Wholesale   &  & 97.7               & 90.9             & 90.9           &  & 96.6               & 89.8             & 90.9           &  & 97.7               & 90.9             & 90.9           \\
Wisconsin   &  & 97.4               & 91.2             & 93.9           &  & 94.7               & 89.5             & 94.7           &  & 93.0               & 95.6             & 95.6           \\ \hline
\end{tabular}}
\caption{\mnote{Predictive performances on the test set: fidelity 
 of the \texttt{MIRET} model
with respect to $\widehat{RF}$, accuracy with respect to the ground truth of the \texttt{MIRET}  and $\widehat{RF}$ models.}}
\label{table: test_acc}
\end{table}

In order to better understand the role of the TE-driven information inserted in the model, we analyze the features' use and the proximity measure among samples, comparing the \texttt{MIRET} tree with respect to the 
$\widehat{RF}$. In particular, in Figure \ref{fig:Feat_selected} we report a one-to-one comparison of each problem and for each $D$ of the level frequency of the features used in the  $\widehat{RF}$  and in the \texttt{MIRET} model. For $\widehat{RF}$, the feature's level frequency is calculated as in \eqref{eq: miret_frequency}.

 For \texttt{MIRET}, the heatmap reports the  feature level frequency as the fraction of times a feature is used in a split at level $d$; thus, it is calculated as
\begin{align*}
f^{\texttt{MIRET}}_{d,j} = \frac{1}{|\Tb(d)|} {\sum_{t\in\B(d)} s^*_{t,j}}.
\end{align*}

 From this Figure, we can also see at a glance the number of features used at each level which gives a rough view of the multi-dimension of the splits of the optimal tree. 
The deeper the $\widehat{RF}$ goes, the more features it tends to use with a low frequency. On the other hand, \texttt{MIRET} uses fewer features at each level.
 Thus, the penalization of the features' use in the objective function allows us to reduce the number of features used in the optimal tree, which are less than those used overall in the forest, and encourages the use of sparse multivariate splits.

\begin{figure}[ht!]
\centering
\includegraphics[width=1.\textwidth]{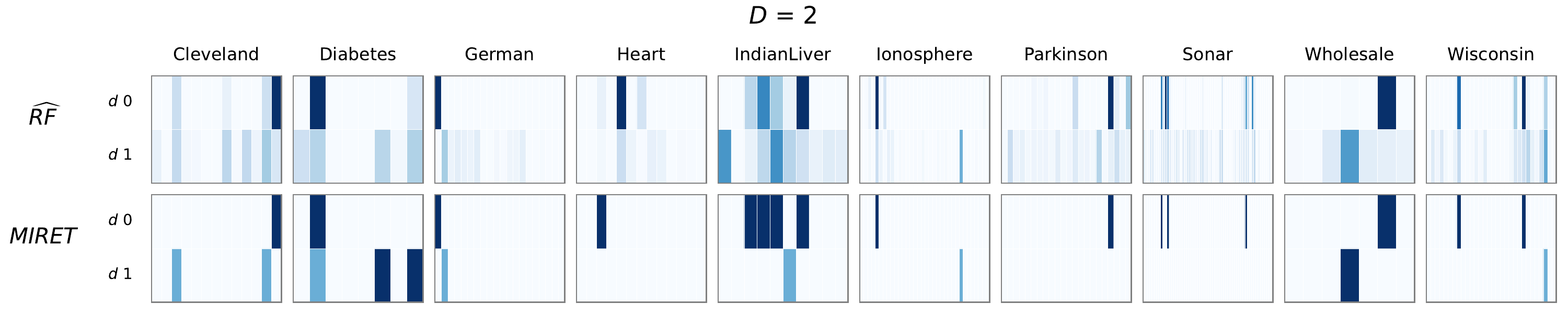}
\vspace{0.1cm}
\includegraphics[width=1.\textwidth]{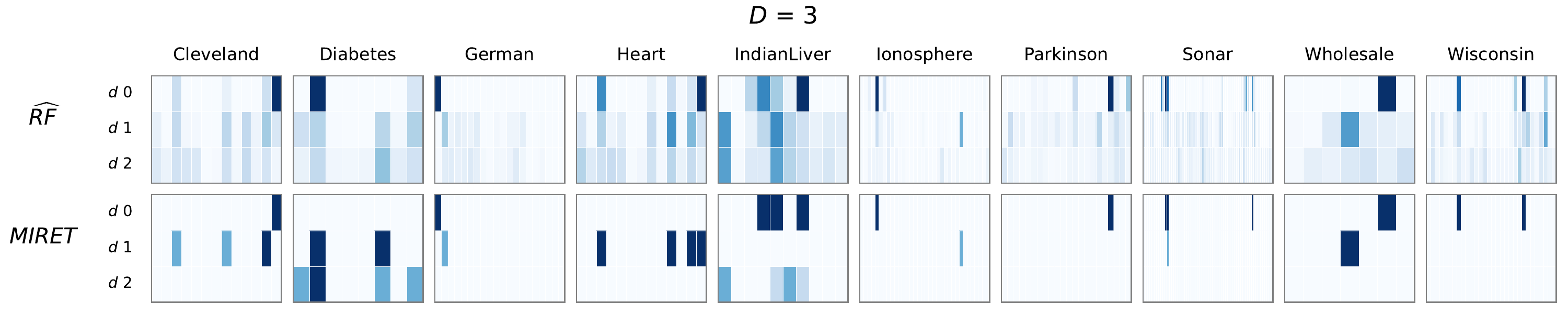}
\vspace{0.1cm}
\includegraphics[width=1.\textwidth]{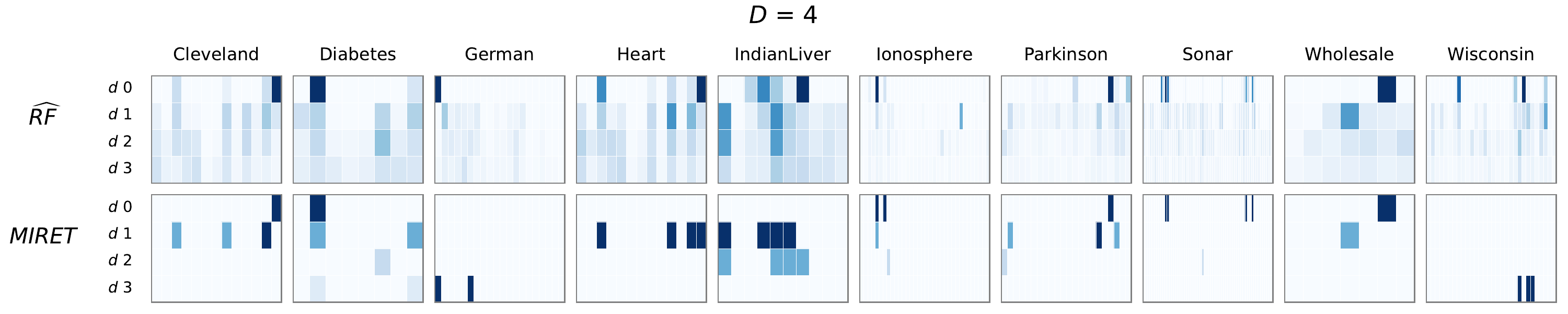} 
\caption{\mnote{Level frequency of features in $\widehat{RF}$ and \texttt{MIRET} model for $D\in\{2,3,4\}$}}
\label{fig:Feat_selected}
\end{figure}
 
As regard the proximity measure, we use pair of samples with a proximity \mnote{ higher than $\overline{m}_{\widehat{RF}}$} in the $\widehat{RF}$ to define hard constraints. Of course, this implies that these samples end up in the same leaves on the training data. This might not be true on the test set. Indeed, the multivariate structure of the splits of \texttt{MIRET}  does not allow for replicating the exact partition of the space made by the $\widehat{RF}$, which uses univariate splits. However, proximity constraints aim to encourage clusters of samples \mnote{that are more likely to end up} in the same leaf in all trees of the $\widehat{RF}$, to belong to the same partition of space. 

\lcomment{Using definition \ref{eq:proximity},  we consider the proximity sets of samples for both the $\widehat{RF}$ and the \texttt{MIRET} tree. Let $\U_{\widehat{RF}}(\overline{m}_{\widehat{RF}})$
 be the set of pairs of samples with proximity measure higher than $\overline{m}_{\widehat{RF}}$
as defined in \eqref{eq:proxiset}, and let $\U_{\texttt{MIRET}}$ be the set of samples defined as:}

\begin{align}
&\U_{\texttt{MIRET}} := \{(i,k)\in\I\times\I: i<k \wedge \text{$z^*_{i,\ell}=z^*_{k,\ell}$ for all $\ell\in\L$} \}.\nonumber
\end{align}

\lnote{Thus, we can calculate the fractions $U_{\widehat{RF}}$,  which measures the percentage of samples that both end up in the same leaf of the target tree ensemble $\widehat{RF}$ at least $\overline{m}_{\widehat{RF}}\%$ of the times and are assigned to the same leaf in \texttt{MIRET}.}
 
\begin{align*}
\lnote{U_{\widehat{RF}}} = \frac{|\U_{\widehat{RF}}(\overline{m}_{\widehat{RF}}) \cap \U_{\texttt{MIRET}}|}{|\U_{\widehat{RF}}(\overline{m}_{\widehat{RF}})|} 
 \end{align*}

\lnote{In Table \ref{table: prox_card_fraction}, we report the cardinality of the  set $\U_{\widehat{RF}}(\overline{m}_{\widehat{RF}})$ both for the training and test set and the value of $U_{\widehat{RF}}$ for the test set. 
The latter metric is not reported for the training set since it is imposed to be maximal with hard constraints in the formulation.
 The higher the values of $U$, the more the \texttt{MIRET} tree mimics the partition of the samples in  $\U_{\widehat{RF}}(\overline{m}_{\widehat{RF}})$  of the $\widehat{RF}$.
 We can observe that pairs of samples from the test set that belong to $\U_{\widehat{RF}}(\overline{m}_{\widehat{RF}})$ consistently end up in the same leaf of the surrogate model. This suggests that the imposed proximity constraints seem to guide the model in imitating the behaviour of the $\widehat{RF}$, even for out-of-sample data.}

\begin{table}[ht!]
\centering
\scalebox{0.7}{
\renewcommand\arraystretch{1.5}
\begin{tabular}{lcccccccccccc}
            &  & \multicolumn{3}{c}{$D = 2$}                                                                                   &  & \multicolumn{3}{c}{$D = 3$}                                                                                   &  & \multicolumn{3}{c}{$D = 4$}                                                                                   \\ \hline
            &  & Train                                              & \multicolumn{2}{c}{Test}                                 &  & Train                                              & \multicolumn{2}{c}{Test}                                 &  & Train                                              & \multicolumn{2}{c}{Test}                                 \\ 
Dataset     &  & $|\U_{\widehat{RF}}(\overline{m}_{\widehat{RF}})|$ & $|\U_{\widehat{RF}}(\overline{m}_{\widehat{RF}})|$ & $U_{\widehat{RF}}$ &  & $|\U_{\widehat{RF}}(\overline{m}_{\widehat{RF}})|$ & $|\U_{\widehat{RF}}(\overline{m}_{\widehat{RF}})|$ & $U_{\widehat{RF}}$ &  & $|\U_{\widehat{RF}}(\overline{m}_{\widehat{RF}})|$ & $|\U_{\widehat{RF}}(\overline{m}_{\widehat{RF}})|$ & $U_{\widehat{RF}}$ \\ \hline
Cleveland   &  & 336                                                & 43                                                 & 100   &  & 69                                                 & 27                                                 & 100   &  & 1014                                               & 55                                                 & 100   \\
Diabetes    &  & 1271                                               & 74                                                 & 100   &  & 315                                                & 3                                                  & 100   &  & 93                                                 & 0                                                  & -   \\
German      &  & 1667                                               & 105                                                & 100   &  & 24                                                 & 3                                                  & 100   &  & 1590                                               & 149                                                & 100   \\
Heart       &  & 1284                                               & {78}                                               & 100   &  & 45                                                 & 100                                                  & 100   &  & 467                                                & 22                                                 & 100   \\
IndianLiver &  & 170                                                & 4                                                  & 100   &  & 2606                                               & 151                                                & 100   &  & 8                                                  & 0                                                  & -   \\
Ionosphere  &  & 6515                                               & 311                                                & 100   &  & 1887                                               & 53                                                 & 100   &  & 8341                                               & 401                                                & 100   \\
Parkinson   &  & 246                                                & 5                                                  & 100   &  & 176                                                & 4                                                  & 100   &  & 1826                                               & 50                                                 & 100   \\
Sonar       &  & 279                                                & {5}                                                & 100   &  & 1                                                  & 0                                                  & -   &  & 0                                                  & 0                                                  & -   \\
Wholesale   &  & 3971                                               & 388                                                & 100   &  & 1743                                               & 112                                                & 100   &  & 12919                                              & 1124                                               & 100   \\
Wisconsin   &  & 26217                                              & 1415                                               & 100   &  & 30427                                            & 1485                                               & 100   &  & 16351                                              & 991                                                & 100   \\ \hline
\end{tabular}}
\caption{\mnote{Cardinality of $\U_{\widehat{RF}}(\overline{m}_{\widehat{RF}})$ for both training and test set and the $U_{\widehat{RF}}$ (in \%) evaluated on the test set. "-" indicates that it is not possible to evaluate $U_{\widehat{RF}}$ since $\U_{\widehat{RF}}(\overline{m}_{\widehat{RF}})$=0. }}
\label{table: prox_card_fraction}
\end{table}

As a final example of the  \texttt{MIRET} model, in Figure 
\ref{fig:tree_cleveland_def} we report the tree generated with $D=3$ on the Cleveland problem, which we visualized using  \texttt{VITE} with the same setting of $\widehat{RF}$
in Figure 
\ref{fig:tree_heatmap}. We observe the actual depth of the tree is 2 although $D=3$. Indeed the nodes at level 2 generate only dummy children that are used only to define the class of all the samples. 
Thus we have the following classification function
$$T^*(x)=\begin{cases}
1  &\text{ if }  (x_{12}\le 0\wedge  -0.009x_2-0.003x_{11}\leq-0.01 ) \vee (x_{12}> 0\ \wedge\  -x_7+x_{11}>-0.641 ) \cr
-1 &\text{ if }  (x_{12}\le 0\wedge  -0.009x_2-0.003x_{11} > -0.01 ) \vee (x_{12}> 0\ \wedge\  -x_7+x_{11}\le -0.641 ) 
\end{cases}
$$

It happens on several problems that the \texttt{MIRET} tree has an actual depth smaller than the maximum depth $D$ of $\widehat{RF}$. This can be easily derived from Figure \ref{fig:Feat_selected}, where the presence of a level $d$ with no "colored" features in the \texttt{MIRET} level frequency plot means that $f^{\texttt{MIRET}}_{d,j} = 0 $ $ \forall j \in \J $, i.e. none of the nodes at level $d$ applies a splitting rule. Thus, the size of the \texttt{MIRET} tree is often smaller than the maximum possible, leading to more interpretable trees. \bigskip

\lnote{
The \texttt{MIRET} methodology is versatile and compatible with Random Forests and boosted ensemble techniques. While we focus largely on the former, we also extend the evaluation to XGBoost models to assess the performances of \texttt{MIRET} when it is applied to identify an optimal representer tree $T^*$ of an XGBoost target model ($\widehat{XGB}$). We report these result in the \ref{sec: miret in xgb}.
}

\begin{figure}
\centering 
\scalebox{0.5}{\tikzset{
  branchnode/.style = {line width=1.8pt, scale=1.3, very thick,          rectangle, rounded corners, draw=dblue, minimum height=1.2cm, text centered, minimum width=1.7cm, font=\normalsize},
  1dummy/.style       = {font=\Large},
  rect/.style = {rectangle, align = left, draw = white, font=\huge, text width=4cm}
}

\tikzset{
    ncbar angle/.initial=90,
    ncbar/.style={
        to path=(\tikztostart)
        -- ($(\tikztostart)!#1!\pgfkeysvalueof{/tikz/ncbar angle}:(\tikztotarget)$)
        -- ($(\tikztotarget)!($(\tikztostart)!#1!\pgfkeysvalueof{/tikz/ncbar angle}:(\tikztotarget)$)!\pgfkeysvalueof{/tikz/ncbar angle}:(\tikztostart)$)
        -- (\tikztotarget)
    },
    ncbar/.default=0.3cm,
}

\tikzset{square left brace/.style={ncbar=0.3cm}}
\tikzset{square right brace/.style={ncbar=-0.3cm}}

\begin{tikzpicture}
  [
    grow                    = down,
    level 1/.style          = {sibling distance=14cm},
    level 2/.style          = {sibling distance=6.5cm},
    level 3/.style          = {sibling distance=3.4cm},
    level distance          = 3.8cm,
    edge from parent/.style = {thick, draw, edge from parent path={(\tikzparentnode) -- (\tikzchildnode)}},
  ]
\node [branchnode] (0) {$[128,109]$}           
    child {
    node[branchnode] (1) {$[101,27]$}   
    child { node [branchnode] (3) {$[3,14]$}    
    child { node [branchnode, draw = mygreen] (7) {$[0,0]$}   
         edge from parent node[1dummy, pos=0.4, left=0.4cm] (e) {$\leq$} }
    child { node [branchnode, draw = mygreen] (8) {$[3,14]$}   
      edge from parent node[1dummy, pos=0.4, right=0.4cm] (k) {$0>-1$}}
    edge from parent node[1dummy, pos=0.4, left=0.4cm] (c) {$-0.009x_2-0.003x_{11}\leq-0.01$} }
    child { node [branchnode] (4) {$[98,13]$}    
    child { node [branchnode, draw = mygreen] (9) {$[98,13]$}   
    edge from parent node[1dummy, pos=0.4, left=0.4cm] (h) {$0\leq1$} }
    child { node [branchnode, draw = mygreen] (10) {$[0,0]$}   
      edge from parent node[1dummy, pos=0.4, right=0.4cm] (k) {$>$}}
    edge from parent node[1dummy, pos=0.4, right=0.4cm] {$>$} }
  edge from parent node [1dummy, pos=0.3, left=0.9cm] (a) {$x_{12}\leq 0$}
}
    child {
     node[branchnode] (2) {$[27,82]$}   
    child { node  [branchnode] (5) {$[14,5]$}    
    child { node [branchnode, draw = mygreen] (11) {$[14,5]$}   
        edge from parent node[1dummy, pos=0.4, left=0.4cm] (j) {$0\leq1$} }
    child { node [branchnode, draw = mygreen] (12) {$[0,0]$}   
        edge from parent node[1dummy, pos=0.4, right=0.4cm] (k) {$>$} }
     edge from parent 
     node [1dummy, pos=0.3, left=0.4cm] {$\leq$}}
    child { node [branchnode] (6) {$[13,77]$}   
    child { node [branchnode, draw = mygreen] (13) {$[0,0]$}   
      edge from parent node[1dummy, pos=0.4, left=0.4cm] (l) {$\leq$}}
    child { node [branchnode, draw = mygreen] (14) {$[13,77]$}   
        edge from parent node [1dummy, pos=0.4, right=0.4cm] (f) {$0>-1$}}
    edge from parent 
     node[1dummy, pos=0.4, right=0.4cm, draw=none] (d) {$-x_7+x_{11}>-0.641$}}
    edge from parent node [1dummy, pos=0.3, right=0.9cm] (b) {$x_{12}> 0$}
};









\node[1dummy, below = 0.03cm of 7]  {class $-1$};
\node[1dummy, below = 0.03cm of 8]  {class $1$};
\node[1dummy, below = 0.03cm of 9]  {class $-1$};
\node[1dummy, below = 0.03cm of 10]  {class $1$};
\node[1dummy, below = 0.03cm of 11]  {class $-1$};
\node[1dummy, below = 0.03cm of 12]  {class $1$};
\node[1dummy, below = 0.03cm of 13]  {class $-1$};
\node[1dummy, below = 0.03cm of 14]  {class $1$};

\end{tikzpicture}}
\caption{Optimal Tree obtained by \texttt{Miret} with $D=3$ on Cleveland. For each node, we report  the number of samples in the class defined by the ground truth: [{\# negative labels}, {\# positive labels}].}

\label{fig:tree_cleveland_def}
\end{figure}
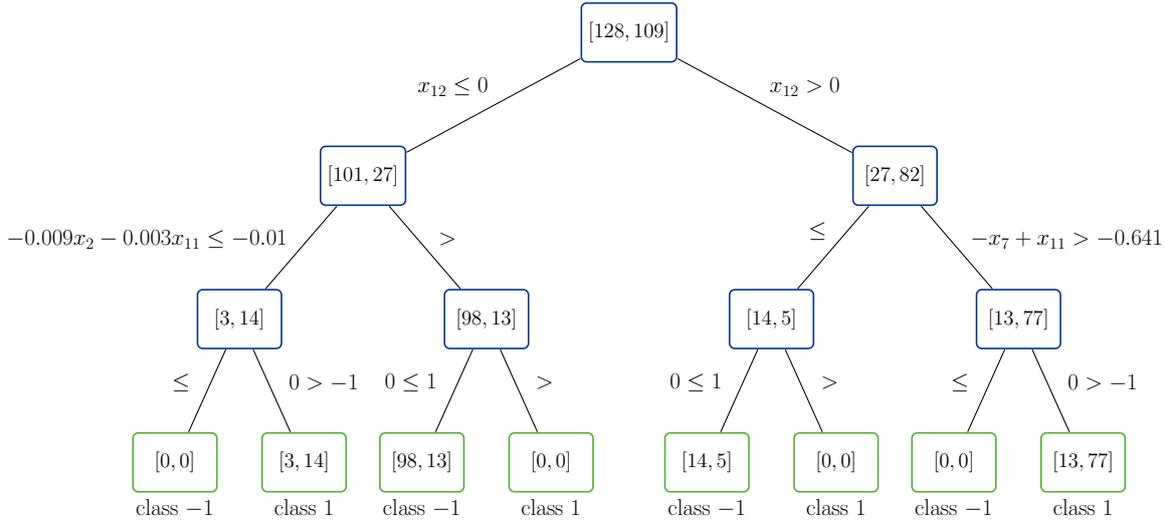


\section{Conclusion}

The paper falls in the field of interpretable representation of a tree-ensemble model, which aims to provide valuable insights into the relationship between the input features and the TE outcomes.
Our contribution is twofold.
Firstly, we propose a visualization tool \texttt{VITE} for ensemble tree models, which allows the user to capture the hierarchical role of features in determining predictions by showing the features' frequency use in the forest. The proposed tool is an addition to the existing visualization tools by helping to understand how features are used in the black-box tree-ensemble model.

Further, we present a mixed-integer linear formulation for learning an interpretable re-built tree (\texttt{MIRET}) from a target tree ensemble model (TE).
\texttt{MIRET} is a multivariate tree with assigned depth $D$, which optimizes a weighted combination of fidelity to the TE and the number of the features used across the tree, gaining in interpretability at a fixed complexity of the tree. In order to improve consistency with the TE, we extracted information from it, such as level frequencies, proximity measures, and class probabilities, and we embed them into the MILP model.
In this paper, we fixed the depth $D$ of the \texttt{MIRET} tree to the one of the target TE. However, in principle, by adapting the definition of level frequencies, we can develop optimal trees with any desired depth.



Results on benchmark datasets show that the proposed model is effective in feature selection and yields a shallow interpretable tree while accurately approximating the tree-ensemble decision function. The \texttt{MIRET} model offers improved interpretability yielding a single optimal tree with intuitive decision paths, which fairly closely replicates the predictive 
capabilities of the target random forest model. 

\section*{Acknowledgements}
\lnote{This research has been partially carried out in the framework of the CADUCEO project (No. F/180025/01-05/X43), supported by the Italian Ministry of Enterprises and Made in Italy. This support is gratefully acknowledged. 
Laura Palagi acknowledges financial support from Progetto di Ricerca Medio Sapienza Uniroma1 - n. RM1221816BAE8A79.
Marta Monaci acknowledges financial support from Progetto Avvio alla Ricerca Sapienza Uniroma1 - n. AR1221816C6DC246.
Giulia Di Teodoro acknowledges financial support from Progetto Avvio alla Ricerca Sapienza Uniroma1 - n. AR12218163024834.}

\bibliographystyle{apalike}
\bibliography{bibliography}    

\appendix
\section{Additional results}\label{app: appendix1}

\lnote{In this section we present additional computational results of \texttt{MIRET}.}

\subsection{Comparison between \texttt{MIRET} and a pure multivariate tree approach}

\lnote{In the following analysis, we test \texttt{MIRET} against a "pure" multivariate tree model that corresponds to \texttt{MIRET} without the TE-driven information (proximity, features frequency, and class probability) introduced in the formulation. This comparison is aimed at evaluating the added value brought by these metrics in guiding the optimal tree formulation and optimization to imitate the tree ensemble.}

\lnote{The benchmark model is denoted as \texttt{MT} (Multivariate optimal Tree) and it is a MILP model to build an optimal tree  trained on the class labels predicted by the $\widehat{TE}$ in a binary classification setting. Its objective function minimizes a trade-off between the misclassification cost and the sum of the features used across the splits in the tree. 
}

\lnote{In this analysis, we compare \texttt{MT} and \texttt{MIRET} applied to the target $\widehat{RF}$. The hyperparameters used are the ones reported in Table \ref{tab: hyperparameters}
and a time limit of $1$ hour was set for all experiments.}


\lnote{The first comparison is on the optimization performances. 
Table \ref{table: mt_miret_timegap} provides running time and MIP Gap values of the two analysed models. It is evident that \texttt{MIRET} consistently outperforms \texttt{MT} in terms of optimization performances. Indeed, at $D=2$, \texttt{MIRET} certifies the optimal solution in 8 out of 10 datasets, whereas \texttt{MT} does so only in one case. Further, \texttt{MIRET} closes the gap in 4 problems at $D=3$ and in 2 at $D=4$, while \texttt{MT} fails to certify the optimal solution at both depths in all of the problems.
These results show some of the advantages brought by the TE-based techniques embedded into MIRET. In particular, the frequency constraints (\ref{cons: freq}, \ref{cons: freq2})  play a crucial role in diminishing the problem size by reducing the feature space at each node. Similarly,  proximity constraints \ref{cons: prox} also contribute to making the problem easier to solve by constraining some binary variables to be equal. In this way, some clusters of points are forced to end up in the same leaf, limiting the set of feasible solutions.
Therefore, TE-based metrics not only make \texttt{MIRET} easier to solve compared to \texttt{MT} but also indicate the model's adaptability and effectiveness in constructing a surrogate tree.} \\
\lnote{Tables \ref{table: mt_miret_pred_2},\ref{table: mt_miret_pred_3},\ref{table: mt_miret_pred_4} show the comparison on fidelity and accuracy performance between \texttt{MT} and \texttt{MIRET} for $D=2, D=3$ and $D=4$, respectively. Focusing on fidelity in the training set, it is possible to see \texttt{MT} beats \texttt{MIRET} in all cases. Despite this, when looking at fidelity in the test set, \texttt{MIRET}'s performances are strictly higher than \texttt{MT}'s ones in 70\%  of cases and equals in two problems. This indicates that \texttt{MT} overfits on training data and does not generalize well on out-of-sample data.  On the other hand, if we focus on the accuracy of the two models with respect to ground truth, the pattern seems to replicate itself, though less sharply. On the training data \texttt{MT} outperforms \texttt{MIRET} in each case, but on the test set \texttt{MIRET} has better accuracy on 11 of the 30 (36.7\%) problems and the same performance on 5 of them (16.7\%). While the accuracy in the test set of \texttt{MT} is better than or equal to that of \texttt{MIRET} in 46.6\% of the cases, and thus there might seem to be no clear winner. However looking at the number of features used in the splitting-rules, it becomes clear that \texttt{MT} does not meet the goal for which \texttt{MIRET} was proposed.
Table \ref{table: mt_miret_features} reports the cardinality of the set $\mathcal{F}$ of different features used overall in \texttt{MT} and in \texttt{MIRET} and the cardinality of the multi-set $\widehat{\mathcal{F}}$ of all features used in the tree, counted as many times as it is used at each node.
 It is immediately noticeable that \texttt{MT} uses a considerably larger number of features than those used by \texttt{MIRET}, losing interpretability. The number of features used then becomes larger the deeper the tree. Indeed, the goal of the \texttt{MIRET} formulation  is to mimic the behavior of the forest with a single tree, making it interpretable. A decision path with multivariate decision rules involving too many features cannot be viewed as truly interpretable. }

\lnote{In Table \ref{table: mt_miret_prox}, we report the cardinalities of the set $\U_{\widehat{RF}}(\overline{m}_{\widehat{RF}})$ on the test set, along with the corresponding $\U$ for both \texttt{MIRET} and \texttt{MT}. These results show that \texttt{MIRET}, differently from \texttt{MT}, consistently assigns pairs of data points with a proximity higher than $\overline{m}_{\widehat{RF}}$ to the same leaf, thereby aligning more closely with $\widehat{RF}$.  
It is important to point out that our primary goal is not only to replicate the $\widehat{RF}$ predictions using a single interpretable tree but also to mirror specific characteristics of $\widehat{RF}$, such as encouraging the use of the most frequent features and approximating, to some extent, the spatial distribution of points returned by the $\widehat{RF}$.}

 
\begin{table}[ht!]
\centering
\scalebox{0.8}{
\renewcommand\arraystretch{1.3}
\begin{tabular}{lcccccccccccccc}
            & \multicolumn{4}{c}{$D=2$}                                  &  & \multicolumn{4}{c}{$D=3$}                                  &  & \multicolumn{4}{c}{$D=4$}                                           \\ \hline
            & \multicolumn{2}{c}{MT} & \multicolumn{2}{c}{MIRET}        &  & \multicolumn{2}{c}{MT} & \multicolumn{2}{c}{MIRET}        &  & \multicolumn{2}{c}{MT}          & \multicolumn{2}{c}{MIRET}        \\ 
Dataset     & Time             & Gap  & Time             & Gap           &  & Time             & Gap  & Time             & Gap           &  & Time             & Gap           & Time             & Gap           \\ \hline
Cleveland   & \underline{3600} & 85.7 & \textbf{44.6}    & \textbf{0.0}  &  & \underline{3600} & 90.9 & \textbf{1942.2}  & \textbf{0.0}  &  & \underline{3600} & 95.0          & \underline{3600} & \textbf{71.5} \\
Diabetes    & \underline{3600} & 99.1 & \textbf{1373.3}  & \textbf{0.0}  &  & \underline{3600} & 98.7 & \underline{3600} & \textbf{96.1} &  & \underline{3600} & \textbf{98.6}          & \underline{3600} & \textbf{98.6} \\
German      & \underline{3600} & 85.7 & \textbf{12.3}    & \textbf{0.0}  &  & \underline{3600} & 98.1 & \underline{3600} & \textbf{73.6} &  & \underline{3600} & 99.2          & \underline{3600} & \textbf{93.4} \\
Heart       & \underline{3600} & 90.0 & \textbf{0.7}     & \textbf{0.0}  &  & \underline{3600} & 90.9 & \underline{3600} & \textbf{73.2} &  & \underline{3600} & 95.7          & \underline{3600} & \textbf{95.0} \\
IndianLiver & \underline{3600} & 95.8 & \underline{3600} & \textbf{75.6} &  & \underline{3600} & 97.8 & \underline{3600} & \textbf{77.1} &  & \underline{3600} & 98.7          & \underline{3600} & \textbf{98.3} \\
Ionosphere  & \underline{3600} & 91.7 & \textbf{3.4}     & \textbf{0.0}  &  & \underline{3600} & 92.9 & \textbf{204.4}   & \textbf{0.0}  &  & \underline{3600} & 94.4          & \textbf{947.6}            & \textbf{0.0}  \\
Parkinson   & \underline{3600} & 83.3 & \textbf{4.6}     & \textbf{0.0}  &  & \underline{3600} & 88.9 & \underline{3600} & \textbf{33.8} &  & \underline{3600} & \textbf{91.7} & \underline{3600} & 96.3          \\
Sonar       & \underline{3600} & 93.7 & \underline{3600} & \textbf{12.0} &  & \underline{3600} & 95.8 & \underline{3600} & \textbf{93.2} &  & \underline{3600} & 97.4          & \underline{3600} & \textbf{97.1} \\
Wholesale   & 284.6            & \textbf{0.0}  & \textbf{0.7}     & \textbf{0.0}           &  & \underline{3600} & 75.0 & \textbf{18.0}    & \textbf{0.0}  &  & \underline{3600} & 94.7          & \textbf{554.2}            & \textbf{0.0}  \\
Wisconsin   & \underline{3600} & 85.7 & \textbf{8.2}     & \textbf{0.0}  &  & \underline{3600} & 88.9 & \textbf{385.5}            & \textbf{0.0}  &  & \underline{3600} & \textbf{85.7} & \underline{3600} & 98.1          \\ \hline
\end{tabular}}
\caption{\mnote{{Comparison {between} \texttt{MT} and \texttt{MIRET} models on
computational times (s) and MIP Gap values (\%) {for RF models}; in boldface the winning values.}}}
\label{table: mt_miret_timegap}
\end{table}

\begin{table}[ht!]
\centering
\scalebox{0.8}{
\renewcommand\arraystretch{1.3}
\begin{tabular}{lccccccccccc}
            & \multicolumn{11}{c}{$D=2$}                                                                                                                             \\ \hline
            & \multicolumn{4}{c}{FID-$\widehat{RF}$}                     &  & \multicolumn{6}{c}{ACCURACY}                                                           \\ \hline
            & \multicolumn{2}{c}{MT}        & \multicolumn{2}{c}{MIRET} &  & \multicolumn{2}{c}{MT}       & \multicolumn{2}{c}{MIRET}     & \multicolumn{2}{c}{$\widehat{RF}$} \\ 
Dataset     & Train          & Test          & Train   & Test            &  & Train         & Test          & Train         & Test          & Train      & Test      \\ \hline
Cleveland   & \textbf{100.0} & 96.7          & 96.6    & \textbf{98.3}   &  & 81.9          & \textbf{83.3} & \textbf{81.9}          & 78.3          & \textbf{81.9}       & 80.0      \\
Diabetes    & \textbf{96.9}  & \textbf{92.9} & 91.9    & 89.6            &  & \textbf{76.2} & 70.8          & 73.5          & \textbf{75.3} & 75.4       & 68.8      \\
German      & \textbf{100.0} & 98.5          & 98.0    & \textbf{99.0}   &  & 66.0          & \textbf{66.0} & \textbf{66.8} & 65.5          & 66.0       & 65.5      \\
Heart       & \textbf{100.0} & \textbf{92.6} & 77.3    & 77.8            &  & \textbf{84.3} & \textbf{81.5} & 74.5          & 77.8          & 84.3       & 85.2      \\
IndianLiver & \textbf{99.4}  & \textbf{96.6} & 97.2    & 94.8            &  & 67.2          & 63.8          & \textbf{67.6} & \textbf{65.5} & 67.8       & 65.5      \\
Ionosphere  & \textbf{100.0} & 91.5          & 98.6    & \textbf{100.0}  &  & \textbf{91.1} & 85.9          & 89.6          & \textbf{94.4} & 91.1       & 94.4      \\
Parkinson   & \textbf{100.0} & 92.3          & 93.6    & \textbf{97.4}   &  & \textbf{88.5} & \textbf{76.9}          & 84.6          & \textbf{76.9}          & 88.5       & 79.5      \\
Sonar       & \textbf{100.0} & 76.2          & 86.8    & \textbf{90.5}   &  & \textbf{90.4} & 66.7          & 79.5          & \textbf{71.4} & 90.4       & 81.0      \\
Wholesale   & \textbf{100.0} & \textbf{98.9} & 99.7    & 97.7            &  & 92.3          & 89.8          & \textbf{92.6} & \textbf{90.9} & 92.3       & 90.9      \\
Wisconsin   & \textbf{100.0} & 95.6          & 99.1    & \textbf{97.4}   &  & \textbf{96.7} & \textbf{93.0} & 95.8          & 91.2          & 96.7       & 93.9      \\ \hline
\end{tabular}}
\caption{\mnote{{Comparison {between} \texttt{MT} and \texttt{MIRET} models with $D=2$ {on} predictive performances on the training and test set. FID-$\widehat{RF}$: Fidelity with respect to $\widehat{RF}$ (\%); \texttt{MIRET} ACCURACY and $\widehat{RF}$ ACCURACY with respect to the ground truth (\%); in boldface the winning values.}}}
\label{table: mt_miret_pred_2}
\end{table}

\begin{table}[ht!]
\centering
\scalebox{0.8}{
\renewcommand\arraystretch{1.3}
\begin{tabular}{lccccccccccc}
            & \multicolumn{11}{c}{$D=3$}                                                                                                                          \\ \hline
            & \multicolumn{4}{c}{FID-$\widehat{RF}$}                     &  & \multicolumn{6}{c}{ACCURACY}                                                        \\ \hline
            & \multicolumn{2}{c}{MT}        & \multicolumn{2}{c}{MIRET} &  & \multicolumn{2}{c}{MT}        & \multicolumn{2}{c}{MIRET} & \multicolumn{2}{c}{$\widehat{RF}$} \\ 
Dataset     & Train          & Test          & Train   & Test            &  & Train          & Test          & Train   & Test            & Train      & Test      \\ \hline
Cleveland   & \textbf{100.0} & 90.0          & 97.9    & \textbf{95.0}   &  & \textbf{87.8}  & \textbf{81.7} & 85.7    & 76.7            & 87.8       & 81.7      \\
Diabetes    & \textbf{98.7}  & 89.6          & 95.4    & \textbf{95.5}   &  & \textbf{78.5}  & 72.7          & 75.9    & \textbf{77.3}   & 79.2       & 76.6      \\
German      & \textbf{97.8}  & \textbf{97.0}          & 95.4    & \textbf{97.0}            &  & \textbf{70.4}  & \textbf{68.5}          & 68.5    & \textbf{68.5}            & 71.6       & 67.5      \\
Heart       & \textbf{100.0} & 85.2          & 95.4    & \textbf{90.7}   &  & \textbf{89.4}  & \textbf{81.5} & 86.6    & 79.6            & 89.4       & 85.2      \\
IndianLiver & \textbf{99.4}  & 93.9          & 95.7    & \textbf{94.0}   &  & \textbf{70.0}  & \textbf{64.7} & 67.2    & \textbf{64.7}            & 70.6       & 65.5      \\
Ionosphere  & \textbf{99.6}  & 90.1          & 95.0    & \textbf{95.8}   &  & \textbf{94.3}  & 88.7          & 89.6    & \textbf{94.4}   & 94.6       & 93.0      \\
Parkinson   & \textbf{100.0} & 89.7          & 89.1    & \textbf{100.0}  &  & \textbf{95.5}  & \textbf{84.6} & 85.9    & 79.5            & 95.5       & 79.5      \\
Sonar       & \textbf{100.0} & 78.6          & 80.1    & \textbf{83.3}   &  & \textbf{100.0} & 78.6          & 80.1    & \textbf{83.3}   & 100.0      & 81.0      \\
Wholesale   & \textbf{100.0} & \textbf{97.7} & 99.7    & 96.6            &  & \textbf{92.9}  & 88.6          & 92.6    & \textbf{89.8}   & 92.9       & 90.9      \\
Wisconsin   & \textbf{100.0} & 93.9          & 96.9    & \textbf{94.7}   &  & \textbf{98.7}  & \textbf{97.4} & 95.6    & 89.5            & 98.7       & 94.7      \\ \hline
\end{tabular}}
\caption{\mnote{{Comparison {between} \texttt{MT} and \texttt{MIRET} models with $D=3$ {on} predictive performances on the training and test set. FID-$\widehat{RF}$: Fidelity with respect to $\widehat{RF}$ (\%); \texttt{MIRET} ACCURACY and $\widehat{RF}$ ACCURACY with respect to the ground truth (\%); in boldface the winning values.}}}
\label{table: mt_miret_pred_3}
\end{table}

\begin{table}[ht!]
\centering
\scalebox{0.8}{
\renewcommand\arraystretch{1.3}
\begin{tabular}{lccccccccccc}
            & \multicolumn{11}{c}{$D=4$}                                                                                                                          \\ \hline
            & \multicolumn{4}{c}{FID-$\widehat{RF}$}                     &  & \multicolumn{6}{c}{ACCURACY}                                                        \\ \hline
            & \multicolumn{2}{c}{MT}        & \multicolumn{2}{c}{MIRET} &  & \multicolumn{2}{c}{MT}        & \multicolumn{2}{c}{MIRET} & \multicolumn{2}{c}{$\widehat{RF}$} \\
Dataset     & Train          & Test          & Train   & Test            &  & Train          & Test          & Train   & Test            & Train      & Test      \\ \hline
Cleveland   & \textbf{100.0} & 88.3          & 93.7    & \textbf{93.3}   &  & \textbf{92.0}  & \textbf{78.3} & 85.7    & 76.7            & 92.0       & 83.3      \\
Diabetes    & \textbf{96.4}  & 91.6          & 92.4    & \textbf{94.2}   &  & \textbf{80.0}  & 74.7          & 76.9    & \textbf{78.6}   & 83.2       & 76.6      \\
German      & \textbf{95.3}  & 85.5          & 87.8    & \textbf{87.0}   &  & \textbf{75.8}  & \textbf{72.5} & 68.8    & 68.0            & 78.0       & 72.0      \\
Heart       & \textbf{99.1}  & 88.9          & 93.1    & \textbf{90.7}   &  & \textbf{93.1}  & \textbf{81.5} & 87.0    & 79.6            & 94.0       & 85.2      \\
IndianLiver & \textbf{98.1}  & 94.8          & 93.5    & \textbf{96.6}   &  & \textbf{75.4}  & \textbf{65.5} & 72.1    & 63.8            & 77.3       & 63.8      \\
Ionosphere  & \textbf{100.0} & \textbf{88.7} & 95.0    & 87.3            &  & \textbf{96.8}  & 87.3          & 92.5    & \textbf{88.7}   & 96.8       & 93.0      \\
Parkinson   & \textbf{100.0} & \textbf{97.4}          & 96.2    & \textbf{97.4}            &  & \textbf{98.7}  & \textbf{84.6}          & 94.9    & \textbf{84.6}            & 98.7       & 82.1      \\
Sonar       & \textbf{100.0} & 78.6          & 81.9    & \textbf{88.1}   &  & \textbf{100.0} & \textbf{76.2} & 81.9    & 71.4            & 100.0      & 78.6      \\
Wholesale   & \textbf{99.4}  & 95.5          & 97.2    & \textbf{97.7}   &  & \textbf{95.2}  & 90.9          & 92.3    & \textbf{90.9}            & 95.2       & \textbf{90.9}      \\
Wisconsin   & \textbf{100.0} & \textbf{93.9}         & 98.5    & 93.0            &  & \textbf{99.3}  & 96.5          & 97.8    & \textbf{95.6}            & 99.3       & \textbf{95.6}      \\ \hline
\end{tabular}}
\caption{\mnote{{Comparison {between} \texttt{MT} and \texttt{MIRET} models with $D=4$ {on} predictive performances on the training and test set. FID-$\widehat{RF}$: Fidelity with respect to $\widehat{RF}$ (\%); \texttt{MIRET} ACCURACY and $\widehat{RF}$ ACCURACY with respect to the ground truth (\%); in boldface the winning values.}}}
\label{table: mt_miret_pred_4}
\end{table}

\begin{table}[ht!]
\centering
\scalebox{0.8}{
\renewcommand\arraystretch{1.3}
\begin{tabular}{lcccccccccccccccc}
            &     &  & \multicolumn{4}{c}{$D=2$}                                                                 &  & \multicolumn{4}{c}{$D=3$}                                                                 &  & \multicolumn{4}{c}{$D=4$}                                                                 \\ \hline
            &     &  & \multicolumn{2}{c}{MT}                      & \multicolumn{2}{c}{MIRET}                   &  & \multicolumn{2}{c}{MT}                      & \multicolumn{2}{c}{MIRET}                   &  & \multicolumn{2}{c}{MT}                      & \multicolumn{2}{c}{MIRET}                   \\
Dataset     & $|\J|$ &  & $|\mathcal{F}|$ & $|\widehat{\mathcal{F}}|$ & $|\mathcal{F}|$ & $|\widehat{\mathcal{F}}|$ &  & $|\mathcal{F}|$ & $|\widehat{\mathcal{F}}|$ & $|\mathcal{F}|$ & $|\widehat{\mathcal{F}}|$ &  & $|\mathcal{F}|$ & $|\widehat{\mathcal{F}}|$ & $|\mathcal{F}|$ & $|\widehat{\mathcal{F}}|$ \\ \hline
Cleveland   & 13  &  & 5               & 7                         & \textbf{3}      & \textbf{3}                &  & 6               & 11                        & \textbf{4}      & \textbf{5}                &  & 11              & 20                        & \textbf{4}      & \textbf{5}                \\
Diabetes    & 8   &  & 7               & 15                        & \textbf{3}      & \textbf{6}                &  & 8               & 37                        & \textbf{4}      & \textbf{7}                &  & 8               & 30                        & \textbf{3}      & \textbf{6}                \\
German      & 20  &  & 7               & 7                         & \textbf{2}      & \textbf{2}                &  & 17              & 29                        & \textbf{2}      & \textbf{2}                &  & 20              & 57                        & \textbf{2}      & \textbf{2}                \\
Heart       & 13  &  & 9               & 10                        & \textbf{1}      & \textbf{1}                &  & 8               & 11                        & \textbf{4}      & \textbf{4}                &  & 11              & 19                        & \textbf{4}      & \textbf{4}                \\
IndianLiver & 10  &  & 7               & 14                        & \textbf{5}      & \textbf{5}                &  & 10              & 31                        & \textbf{5}      & \textbf{9}                &  & 10              & 32                        & \textbf{5}      & \textbf{8}                \\
Ionosphere  & 34  &  & 10              & 12                        & \textbf{2}      & \textbf{2}                &  & 12              & 12                        & \textbf{2}      & \textbf{2}                &  & 15              & 18                        & \textbf{3}      & \textbf{4}                \\
Parkinson   & 22  &  & 6               & 6                         & \textbf{1}      & \textbf{1}                &  & 9               & 9                         & \textbf{1}      & \textbf{1}                &  & 10              & {12}                      & \textbf{5}      & \textbf{6}                \\
Sonar       & 60  &  & 16              & 16                        & \textbf{3}      & \textbf{3}                &  & 20              & 24                        & \textbf{3}      & \textbf{4}                &  & 32              & 38                        & \textbf{5}      & \textbf{5}                \\
Wholesale   & 7   &  & 3               & 4                         & \textbf{2}      & \textbf{3}                &  & 3               & 4                         & \textbf{2}      & \textbf{3}                &  & 7               & 14                        & \textbf{2}      & \textbf{2}                \\
Wisconsin   & 30  &  & 7               & 7                         & \textbf{3}      & \textbf{3}                &  & 9               & 9                         & \textbf{2}      & \textbf{2}                &  & 7               & {7}                       & \textbf{3}      & \textbf{3}                \\ \hline
\end{tabular}}
\caption{\mnote{{Comparison {between} \texttt{MT} and \texttt{MIRET} models {on} the number of features selected. $\mathcal{F}$ is the set of unique features used in the tree and $\widehat{\mathcal{F}}$ is the multi-set of features used in the tree. A feature is present in the multiset as many times as it is used at each node; in boldface the winning values.}}}
\label{table: mt_miret_features}
\end{table}

\begin{table}[ht!]
\centering
\scalebox{0.8}{
\renewcommand\arraystretch{1.3}
\begin{tabular}{lcccccccccccc}
            &  & \multicolumn{3}{c}{$D=2$}                                                                   &  & \multicolumn{3}{c}{$D=3$}                                                                   &  & \multicolumn{3}{c}{$D=4$}                                                                   \\ \hline
            &  & $|\U_{\widehat{RF}}(\overline{m}_{\widehat{RF}})|$ & \multicolumn{2}{c}{$U_{\widehat{RF}}$} &  & $|\U_{\widehat{RF}}(\overline{m}_{\widehat{RF}})|$ & \multicolumn{2}{c}{$U_{\widehat{RF}}$} &  & $|\U_{\widehat{RF}}(\overline{m}_{\widehat{RF}})|$ & \multicolumn{2}{c}{$U_{\widehat{RF}}$} \\
Dataset     &  &                                                    & \texttt{MT}        & \texttt{MIRET}    &  &                                                    & \texttt{MT}        & \texttt{MIRET}    &  &                                                    & \texttt{MT}        & \texttt{MIRET}    \\ \hline
Cleveland   &  & 43                                                 & \textbf{100.0}     & \textbf{100.0}    &  & 27                                                 & \textbf{100.0}     & 77.8              &  & 55                                                 & \textbf{100.0}     & 80.0              \\
Diabetes    &  & 74                                                 & \textbf{100.0}     & 64.9              &  & 3                                                  & \textbf{100.0}     & \textbf{100.0}    &  & 0                                                  & -                  & -                 \\
German      &  & 105                                                & \textbf{100.0}     & 99.0              &  & 3                                                  & \textbf{100.0}     & \textbf{100.0}    &  & 149                                                & \textbf{100.0}     & 67.8              \\
Heart       &  & 78                                                 & \textbf{100.0}     & 65.4              &  & 1                                                  & \textbf{100.0}     & \textbf{100.0}    &  & 22                                                 & \textbf{100.0}     & \textbf{100.0}             \\
IndianLiver &  & 4                                                  & \textbf{100.0}     & \textbf{100.0}    &  & 151                                                & \textbf{100.0}     & 96.8              &  & 0                                                  & -                  & -                 \\
Ionosphere  &  & 311                                                & \textbf{100.0}     & 98.1              &  & 53                                                 & \textbf{100.0}     & 83.0              &  & 401                                                & \textbf{100.0}     & 66.8              \\
Parkinson   &  & 5                                                  & \textbf{100.0}     & \textbf{100.0}    &  & 4                                                  & \textbf{100.0}     & \textbf{100.0}    &  & 50                                                 & \textbf{100.0}     & \textbf{100.0}             \\
Sonar       &  & 5                                                  & \textbf{100.0}     & 40.0              &  & 0                                                  & -                  & -                 &  & 0                                                  & -                  & -                 \\
Wholesale   &  & 388                                                & \textbf{100.0}     & \textbf{100.0}    &  & 112                                                & \textbf{100.0}     & \textbf{100.0}    &  & 1124                                               & \textbf{100.0}     & 96.7              \\
Wisconsin   &  & 1415                                               & \textbf{100.0}     & \textbf{100.0}    &  & 1485                                               & \textbf{100.0}     & 96.4              &  & 991                                                & \textbf{100.0}     & 96.0              \\ \cline{1-1} \cline{4-6} \cline{8-10} \hline
\end{tabular}}
\caption{\mnote{{Cardinality of $\U_{\widehat{RF}}(\overline{m}_{\widehat{RF}})$ for the test set and comparison {between} \texttt{MT} and \texttt{MIRET} models {on} the $U_{\widehat{RF}}$ (in \%) evaluated on the test set. - indicates that it is not possible to evaluate $U_{\widehat{RF}}$ since $\U_{\widehat{RF}}(\overline{m}_{\widehat{RF}})$=0.; in boldface the winning values.}}}
\label{table: mt_miret_prox}
\end{table}

\subsection{\lnote{\texttt{MIRET} applied to XGBoost target models.}} \label{sec: miret in xgb}
\lnote{In the following analysis, we assess \texttt{MIRET} performance when applied to gradient boosting methods, using a target XGBoost model ($\widehat{XGB}$).
In particular, we use the XGBoost method as implemented in \texttt{xgboost.XGBClassifier} with the following settings:
\begin{itemize}
    \item maximum depth $D\in\{2,3,4\}$;
    \item $|\E| = 100$;
    \item weights $w^e, e\in \E$ as defined in \cite{Tan2020}.
\end{itemize}
To select the set of hyperparameters, we performed a grid-search within a $4$-fold cross-validation for $\overline{m}_{TE}$, $\gamma_d$, $d\in\D$, and $\alpha$, as in the case of Random Forest indicated in subsection \ref{sub: hyperparam_set}.\\
Table \ref{table: gap time miret xgboost} presents the computational time and MIP Gap values results. Similarly to the results obtained applying \texttt{MIRET} to $\widehat{RF}$, we can certify optimality for all but one problem at $D=2$, and for the majority of problems at $D=3$. At $D=4$, optimality is certified for three out of the ten problems examined. This is still a noteworthy achievement, given that the scalability of optimal trees diminishes as their maximum depth increases.
Further, Tables \ref{table: train_acc_xgb} and \ref{table: test_acc_xgb} reveal that \texttt{MIRET} maintains a high level of fidelity to the $\widehat{XGB}$ models across various depths, with an average fidelity score for the training set of 90.0\% at $D=2$, 91.6\% at $D=3$ and 89.2\% at $D=4$ and with an average fidelity score for the test set of 89.3\% at a $D=2$, 91.4\% at $D=3$ and 89.6\% at $D=4$. Regarding performance metrics, it is evident that MIRET's accuracy mirrors that of $\widehat{XGB}$ models quite well over all datasets, particularly in the test set. This suggests that our model retains the generalization capabilities of the target $\widehat{XGB}$ and offers the added advantage of improved interpretability.
Table \ref{table: prox_card_fraction_xgb} reports the metrics based on proximity measure, as described in Section \ref{sec: results}. We can observe that pairs of samples from the test set in $\U_{\widehat{XGB}}(\overline{m}_{\widehat{XGB}})$ are also assigned to the same leaf of the surrogate model in most of the cases. The high values of $U_{\widehat{XGB}}$ suggest that the proximity constraints imposed in \texttt{MIRET} seem to guide the model in imitating the behavior of the $\widehat{RF}$, even for out-of-sample data.}

\begin{table}[ht!]
\centering
\scalebox{1}{
\renewcommand\arraystretch{1.3}
\begin{tabular}{lccccccccc}
            & \multicolumn{1}{l}{} & \multicolumn{2}{c}{$D=2$}     & \multicolumn{1}{l}{} & \multicolumn{2}{c}{$D=3$}     & \multicolumn{1}{l}{} & \multicolumn{2}{c}{$D=4$}     \\\hline
Dataset     &                      & Time             & Gap        &                      & Time             & Gap        &                      & Time             & Gap        \\ \hline
Cleveland   &                      & 2.3              & \textbf{0} &                      & 36.6             & \textbf{0} &                      & \underline{3600} & 92.6       \\ 
Diabetes    &                      & 0.8              & \textbf{0} &                      & \underline{3600} & 92.8       &                      & \underline{3600} & 99.5       \\
German      &                      & \underline{3600} & 96.1       &                      & \underline{3600} & 98.9       &                      & \underline{3600} & 98.0       \\
Heart       &                      & 168.6            & \textbf{0} &                      & 138.6            & \textbf{0} &                      & \underline{3600} & 77.2       \\
IndianLiver &                      & 55.8             & \textbf{0} &                      & 27.00            & \textbf{0} &                      & \underline{3600} & 98.3       \\
Ionosphere  &                      & 21.3             & \textbf{0} &                      & 151.2            & \textbf{0} &                      & 714.5            & \textbf{0} \\
Parkinson   &                      & 21.5             & \textbf{0} &                      & 124.3            & \textbf{0} &                      & 543.6            & \textbf{0} \\
Sonar       &                      & 108.4            & \textbf{0} &                      & \underline{3600} & 96.9       &                      & \underline{3600} & 98.7       \\
Wholesale   &                      & 0.9              & \textbf{0} &                      & 4.4              & \textbf{0} &                      & \underline{3600} & 87.1       \\
Wisconsin   &                      & 7.4              & \textbf{0} &                      & \underline{3600} & 28.6       &                      & 245.0            & \textbf{0} \\ \hline
\end{tabular}}
\caption{\gnote{{\texttt{MIRET} models in
computational times (s) and MIP Gap values (\%) for XGBoost models.}}}
\label{table: gap time miret xgboost}
\end{table}

\begin{table}[ht!]
\centering
\scalebox{0.85}{
\renewcommand\arraystretch{1.3}
\begin{tabular}{lcccccccccccc}
            &  & \multicolumn{3}{c}{$D = 2$}                            &  & \multicolumn{3}{c}{$D = 3$}                            &  & \multicolumn{3}{c}{$D = 4$}                            \\ \hline
            &  &                    & \multicolumn{2}{c}{ACCURACY}      &  &                    & \multicolumn{2}{c}{ACCURACY}      &  &                    & \multicolumn{2}{c}{ACCURACY}      \\
Dataset     &  & FID-$\widehat{XGB}$ & $\texttt{MIRET}$ & $\widehat{XGB}$ &  & FID-$\widehat{XGB}$ & $\texttt{MIRET}$ & $\widehat{XGB}$ &  & FID-$\widehat{XGB}$ & $\texttt{MIRET}$ & $\widehat{XGB}$ \\ \hline
Cleveland   &  & 94.5              & 85.2             & 89.0           &  & 97.0               & 86.1             & 88.2           &  & 87.8               & 84.8             & 95.4           \\
Diabetes    &  & 89.1               & 76.1             & 80.8           &  & 87.9               & 77.4             & 82.2           &  & 88.8               & 84.8             & 92.7           \\
German      &  & 84.9               & 76.5             & 86.1           &  & 86.1               & 77.6             & 89.8           &  & 78.8               & 74.3             & 94.3           \\
Heart       &  & 93.5               & 85.6             & 88.4           &  & 90.3               & 83.8             & 87.0           &  & 87.5               & 83.8             & 92.6           \\
IndianLiver &  & 82.1               & 74.9             & 88.1           &  & 90.5               & 74.7             & 79.0           &  & 83.8               & 76.5             & 89.6           \\
Ionosphere  &  & 94.3              & 93.2             & 98.9           &  & 95.4               & 95.4             & 100.0           &  & 92.1               & 91.8             & 99.6           \\
Parkinson   &  & 90.4               & 90.4             & 98.7           &  & 96.2              & 96.2             & 100.0           &  & 92.3               & 92.3             & 100           \\
Sonar       &  & 78.3               & 78.3             & 100.0           &  & 83.1               & 82.5             & 99.4           &  & 90.4               & 90.4            & 100.0           \\
Wholesale   &  & 94.3               & 91.5             & 97.2           &  & 94.6               & 94.0             & 99.4           &  & 93.5               & 93.5             & 100.0           \\
Wisconsin   &  & 98.5               & 97.8             & 99.3           &  & 95.2               & 95.2             & 100.0           &  & 96.7               & 96.5             & 99.8           \\ \hline
\end{tabular}}
\caption{\lcomment{Predictive performances on the training set:}
\lnote{ fidelity  of the \texttt{MIRET} model
with respect to $\widehat{XGB}$, accuracy with respect to the ground truth of the \texttt{MIRET}  and $\widehat{XGB}$ models.}}
\label{table: train_acc_xgb}
\end{table}

\begin{table}
\centering
\scalebox{0.85}{
\renewcommand\arraystretch{1.3}
\begin{tabular}{lcccccccccccc}
            &  & \multicolumn{3}{c}{$D = 2$}                            &  & \multicolumn{3}{c}{$D = 3$}                            &  & \multicolumn{3}{c}{$D = 4$}                            \\ \hline
            &  &                    & \multicolumn{2}{c}{ACCURACY}      &  &                    & \multicolumn{2}{c}{ACCURACY}      &  &                    & \multicolumn{2}{c}{ACCURACY}      \\
Dataset     &  & FID-$\widehat{XGB}$ & $\texttt{MIRET}$ & $\widehat{XGB}$ &  & FID-$\widehat{XGB}$ & $\texttt{MIRET}$ & $\widehat{XGB}$ &  & FID-$\widehat{XGB}$ & $\texttt{MIRET}$ & $\widehat{XGB}$ \\ \hline
Cleveland   &  & 90.0              & 78.3             & 88.3           &  & 90.0               & 81.7             & 85.0           &  & 91.7               & 85.0             & 86.7           \\
Diabetes    &  & 90.3               & 76.6             & 81.2           &  & 91.6               & 77.9             & 78.6           &  & 83.1               & 73.4             & 76.0           \\
German      &  & 80.5               & 73.0             & 74.5           &  & 85.0               & 71.5             & 76.5           &  & 84.5               & 72.5             & 75.0           \\
Heart       &  & 92.6              & 87.0             & 90.7           &  & 94.4               & 88.9             & 90.7           &  & 94.4               & 88.9             & 90.7           \\
IndianLiver &  & 84.5               & 67.2             & 74.1           &  & 91.4               & 69.0             & 69.0           &  & 85.3               & 69.0             & 68.1           \\
Ionosphere  &  & 91.5              & 88.7             & 85.9           &  & 93.0               & 90.1             & 88.7           &  & 88.7               & 84.5             & 93.0           \\
Parkinson   &  & 92.3               & 84.6             & 87.2           &  & 94.9              & 92.3             & 92.3           &  & 94.9               & 87.2             & 92.3           \\
Sonar       &  & 78.6               & 76.2             & 92.9           &  & 81.0               & 81.0             & 90.5           &  & 78.6               & 76.2            & 88.1           \\
Wholesale   &  & 95.5               & 90.9             & 90.9           &  & 95.5               & 92.0             & 89.8           &  & 97.7               & 88.6             & 90.9           \\
Wisconsin   &  & 93.9               & 95.6             & 98.2           &  & 97.4               & 94.7             & 97.4           &  & 97.4               & 99.1             & 98.2           \\ \hline
\end{tabular}}

\caption{\lcomment{Predictive performances on the test set:}
\lnote{fidelity 
 of the \texttt{MIRET} model
with respect to $\widehat{XGB}$, accuracy with respect to the ground truth of the \texttt{MIRET}  and $\widehat{XGB}$ models.}}
\label{table: test_acc_xgb}
\end{table}

\begin{table}[ht!]
\centering
\scalebox{0.67}{
\renewcommand\arraystretch{1.5}
\begin{tabular}{lcccccccccccc}
            &  & \multicolumn{3}{c}{$D = 2$}                                                                                   &  & \multicolumn{3}{c}{$D = 3$}                                                                                   &  & \multicolumn{3}{c}{$D = 4$}                                                                                   \\ \hline
            &  & Train                                              & \multicolumn{2}{c}{Test}                                 &  & Train                                              & \multicolumn{2}{c}{Test}                                 &  & Train                                              & \multicolumn{2}{c}{Test}                                 \\ 
Dataset     &  & $|\U_{\widehat{XGB}}(\overline{m}_{\widehat{XGB}})|$ & $|\U_{\widehat{XGB}}(\overline{m}_{\widehat{XGB}})|$ & $U_{\widehat{XGB}}$ &  & $|\U_{\widehat{XGB}}(\overline{m}_{\widehat{XGB}})|$ & $|\U_{\widehat{XGB}}(\overline{m}_{\widehat{XGB}})|$ & $U_{\widehat{XGB}}$ &  & $|\U_{\widehat{XGB}}(\overline{m}_{\widehat{XGB}})|$ & $|\U_{\widehat{XGB}}(\overline{m}_{\widehat{XGB}})|$ & $U_{\widehat{XGB}}$ \\ \hline
Cleveland   &  & 1093                                                & 54                                                 & 100   &  & 449                                                 & 31                                                 & 100   &  & 31                                                  & 9                                                  &  100   \\
Diabetes    &  & 7062                                               & 495                                                 & 100   &  & 398                                                & 18                                                  & 94.4   &  & 68                                                  & 2                                                  & 100    \\
German      &  & 166                                              & 14                                                & 100   &  & 17                                                 & -                                                  & -   &  & 23                                                  & 0                                                  &   -  \\
Heart       &  & 321                                               & 15                                               & 100   &  & 200                                                 & 9                                                  & 100   &  & 155                                                  & 5                                                  & 100    \\
IndianLiver &  & 741                                                & 62                                                  & 96.8   &  & 2088                                               & 157                                                & 100   &  & 70                                                  & 30                                                  & 100    \\
Ionosphere  &  & 3328                                               & 80                                                & 98.8   &  & 10255                                               & 248                                                 & 99.6   &  & 5108                                                 & 154                                                 &  100   \\
Parkinson   &  & 1137                                                & 38                                                  & 100   &  & 1810                                                & 68                                                  & 100   &  & 1475                                                  & 50                                                  &  100   \\
Sonar       &  & 117                                                & 2                                                & 100   &  & 31                                                  & 1                                                  &   100  &  & 63                                                  & 0                                                  &  -   \\
Wholesale   &  & 17606                                               & 1262                                                & 100   &  & 13657                                               & 1108                                                & 100   &  & 6001                                                 & 713                                                 & 100    \\
Wisconsin   &  & 31491                                              & 2048                                               & 99.9   &  & 26855                                              & 1630                                               & 100   &  & 20.964                                                & 1702                                                &  100   \\ \hline
\end{tabular}}
\caption{\lnote{Cardinality of $\U_{\widehat{XGB}}(\overline{m}_{\widehat{XGB}})$ for both training and test set and the $U_{\widehat{XGB}}$ (in \%) evaluated on the test set. "-" indicates that it is not possible to evaluate $U_{\widehat{XGB}}$ since $\U_{\widehat{XGB}}(\overline{m}_{\widehat{XGB}})$=0. }}
\label{table: prox_card_fraction_xgb}
\end{table}

\clearpage

\section{Additional tables}\label{app: appendix2}

\text{In this section we present additional tables.} Table \ref{tab: notation} presents a summary of all the notation of sets, parameters, and hyperparameters adopted in the paper, and Tables  \ref{tab: hyperparameters_RF} and  \ref{tab: hyperparameters} report hyperparameters selected with the 4-fold cross validation \lnote{for the $\widehat{RF}$ and $\widehat{XGB}$ target models respectively.}

\begin{table}[ht!]
\centering
\scalebox{1}{
\renewcommand\arraystretch{1.3}
\begin{tabular}{lcccccccccccc}
            &  & \multicolumn{3}{c}{$D = 2$}                              &  & \multicolumn{3}{c}{$D = 3$}                              &  & \multicolumn{3}{c}{$D = 4$}                                 \\ \hline
Dataset     &  & $h$        & $\alpha$    & $\overline{m}_{\widehat{RF}}$ &  & $h$        & $\alpha$    & $\overline{m}_{\widehat{RF}}$ &  & $h$           & $\alpha$    & $\overline{m}_{\widehat{RF}}$ \\ \hline
Cleveland   &  & 100/3      & 0.2         & 1.00                          &  & 100/3      & 0.5         & 1.00                          &  & 100/3 & 0.5 & 0.85                  \\
Diabetes    &  & 100/2      & 0.2         & 1.00                          &  & 100/3      & 0.2         & 1.00                          &  & 100/3         & 0.5         & 1.00                          \\
German      &  & 100/2      & 0.2         & 1.00                          &  & 100/4      & 0.8         & 1.00                          &  & 100/4 & 0.5 & 0.85                 \\
Heart       &  & $0^*$ & 0.2 & 0.85                  &  & 100/2      & 0.2         & 1.00                          &  & 100/3 & 0.5 & 0.85                  \\
IndianLiver &  & $0^*$      & 0.2         & 1.00                          &  & 100/2 & 0.2 & 0.90                 &  & 100/2         & 0.2         & 1.00                          \\
Ionosphere  &  & 100/3      & 0.2         & 1.00                          &  & $0^*$      & 0.5         & 1.00                          &  & 100/2 & 0.4 & 0.85                 \\
Parkinson   &  & 100/4      & 0.5         & 1.00                          &  & 100/3      & 0.5         & 1.00                          &  & $0^*$ & 0.2 & 0.85                 \\
Sonar       &  & $0^*$ & 0.2 & 0.85                  &  & 100/4      & 0.5         & 1.00                          &  & $0^*$ & {0.2} & {1.00}                  \\
Wholesale   &  & 100/4      & 0.2         & 1.00                          &  & 100/2      & 0.2         & 1.00                          &  & {$0^*$} & {0.4} & {0.85}                  \\
Wisconsin   &  & 100/3      & 0.2         & 1.00                          &  & 100/3      & 0.5         & 0.90                          &  & 100/2         & 0.2         & 1.00                          \\ \hline
\end{tabular}}
\caption{\mnote{Hyperparameters for \texttt{MIRET} on $\widehat{RF}$ selected with the 4-fold cross-validation. The apex $^*$ indicates that we set $\gamma_d = 0$.}}
\label{tab: hyperparameters_RF}
\end{table}

\begin{table}[ht!]
\centering
\scalebox{1}{
\renewcommand\arraystretch{1.3}
\begin{tabular}{lcccccccccccc}
            &  & \multicolumn{3}{c}{$D = 2$}                              &  & \multicolumn{3}{c}{$D = 3$}                              &  & \multicolumn{3}{c}{$D = 4$}                                 \\ \hline
Dataset     &  & $h$        & $\alpha$    & $\overline{m}_{\widehat{XGB}}$ &  & $h$        & $\alpha$    & $\overline{m}_{\widehat{XGB}}$ &  & $h$           & $\alpha$    & $\overline{m}_{\widehat{XGB}}$ \\ \hline
Cleveland   &  &  $0^*$      & 0.1         & 0.85                          &  & 100/5      & 0.2         & 0.90                          &  & 100/3 & 0.6 & 0.85                  \\
Diabetes    &  & 100/3      & 0.5         & 0.85                          &  & 100/2      & 0.5         & 0.85                          &  & 100/3         & 0.2         & 0.90                          \\
German      &  &  $0^*$      & 0.2         & 0.90                          &  &  $0^*$     & 0.2         & 0.90                          &  & 100/20 & 0.5 & 0.85                 \\
Heart       &  & 100/3 & 0.2 & 0.85                  &  & 100/2      & 0.5         & 0.85                          &  & 100/2 & 0.5 & 0.85                  \\
IndianLiver &  & 100/3    & 0.4         & 0.85                         &  & 100/2 & 0.5 & 0.85                 &  & 0         & 0.1         & 0.85                          \\
Ionosphere  &  & 100/3      & 0.2         & 0.85                          &  & 100/3      & 0.2         & 0.85                          &  & 100/3 & 0.6 & 0.85                 \\
Parkinson   &  & 100/3      & 0.2         & 0.85                          &  & 100/3      & 0.2         & 0.85                          &  & 100/3 & 0.6 & 0.85                 \\
Sonar       &  & $0^*$ & 0.8 & 0.85                  &  & 100/3      & 0.2         & 0.90                          &  & 100/3 & 0.1 & 0.85                  \\
Wholesale   &  & 100/3      & 0.2         & 0.85                         &  & 100/2      & 0.2         & 0.85                          &  & 100/3 & {0.2} & {0.85}                  \\
Wisconsin   &  & 100/3      & 0.2         & 0.85                         &  & 100/3      & 0.6         & 0.85                          &  & 100/3         & 0.4         & 0.85                          \\ \hline
\end{tabular}}
\caption{\mnote{Hyperparameters for \texttt{MIRET} on $\widehat{XGB}$ selected with the 4-fold cross-validation. The apex $^*$ indicates that we set $\gamma_d = 0$.}}
\label{tab: hyperparameters}
\end{table}

\begin{table}[ht]
\centering
\scalebox{1.}{
\renewcommand\arraystretch{1.5}
\begin{tabular}{lll}
Notation                                    & \multicolumn{1}{c}{} & Description                                                                      \\ \hline\hline
\textbf{Sets}                               &                      &                                                                                  \\
$\D$                                      &                      & Tree levels        \\
$\E$                                    &    & Tree estimators of the $TE$\\
$\Tb$                                       &                      & Branch nodes                         \\
$\Tl$                                       &                      & Leaf nodes                  \\

$\Tb(d)$                                      &                      & Branch nodes at level {$d$} of the tree                                              \\
$\Tbf$                                    &                      & Branch nodes not adjacent to leaves                                              \\
$\Tbl$                                    &                      & Branch nodes adjacent to leaves                                              \\

$\S(t)$                     &                      & Sub-leaves, i.e. leaf nodes of the subtree rooted at node $t\in\Tb$                           \\
$\S_L(t)$    &                      & {Left sub-leaves}, i.e. leaves under the left branch of $t \in \Tb$                                     \\
$\S_R(t)$  &                      & {Right sub-leaves}, i.e. leaves under the right branch of $t \in \Tb$                                 \\
$\I$                                        &                      & Index set of data samples                                                        \\
$\I_t$                                      &                      & Index set of data samples assigned to node $t\in\Tb$                                   \\
$\J$                                        &                      & Index set of features                                                    \\
$\J_{\gamma}(d)$                                        &                      & Index set of features with a level frequency at $d$ in ${TE}$ greater than $\gamma_d$                                                   \\
$\mathcal{U}(\overline{m}_{TE})$                                        &                      & Pairs of samples with a proximity measure in ${TE}$ greater than $\overline{m}_{TE}$                                                   \\
\textbf{Parameters}  &                      &                                                                                  \\
${p}^i$                                         &                      &  Class probability of sample $i$ predicted by ${TE}$                                                     \\
$f_{d,j}$                                         &                      &   Frequency of feature $j$ at level $d$ in ${TE}$                                                           \\
$c_{\ell}$                                         &                      & Class label pre-assigned to leaf node $\ell$                                                               \\
$\varepsilon$                                         &                      & Parameter to model the closed inequality in routing constraints                                                               \\
${m_{i,k}}$                                         &                      & Proximity measure between the pair of samples ($x^i,x^k$)                                                               \\
\textbf{Hyperparameters}                    &                      &                                                                                  \\
$D$                                         &                      & Maximal depth  of the tree                                                                  \\
$\alpha$                     &                      & Penalty parameter for the feature selection       \\
$\gamma_d$                     &                      & Frequency threshold used to determine $\J_{\gamma}(d)$       \\
$\overline{m}_{TE}$                     &                      & Proximity threshold used to determine $\U(\overline{m}_{TE})$       \\
   \hline         
\end{tabular}}
\caption{Notation: sets, parameters and hyperparameters.}
\label{tab: notation}
\end{table}

\end{document}